\documentclass[11pt,a4paper, reqno]{amsart}

\usepackage{latexsym}

\usepackage{amsmath,amsfonts,amssymb,amsthm,amscd,epsfig}
\usepackage[english]{babel}
\usepackage{indentfirst}
\usepackage[mathscr]{eucal}
\usepackage{graphics}

\numberwithin{equation}{section}

\usepackage{amsfonts,layout}
\newcommand{\Z}{\mathbb Z}
\newcommand{\Q}{\mathbb Q}
\newcommand{\R}{\mathbb R}
\def\QED{\begin{flushright}$\Box$\end{flushright}}

\def \lg{\langle}
\def \rg{\rangle}
\newcommand{\be}{\begin{equation}}
\newcommand{\ee}{\end{equation}}

\newcommand{\commentout}[1]{}
\newcommand{\norm}[1]{\lVert #1 \rVert}
\newcommand{\abs}[1]{\lvert #1 \rvert}
\newcommand{\br}{\begin{eqnarray}}
\newcommand{\er}{\end{eqnarray}}

\def\ep{\varepsilon}

\def\ds{\displaystyle}

\newtheorem{Theorem}{Theorem}[section]

\newtheorem{Lemma}[Theorem]{Lemma}
\newtheorem{Corollary}[Theorem]{Corollary}

\title[Homogenization and enhancement for the $G-$equation]
{Homogenization and enhancement for the $G-$equation}
\author[P. Cardaliaguet, J. Nolen and P. E. Souganidis]
{P. Cardaliaguet \and J. Nolen \and P. E. Souganidis}

\address{D\'epartement de Math\'ematiques, Universit\'e de Brest,
29285 BREST Cedex,
FRANCE}
\email{Pierre.Cardaliaguet@univ-brest.fr }
\address{Department of Mathematics, Duke University, Durham, NC 27708, USA}
\email{nolen@math.duke.edu}
\address{Department of Mathematics, University of Chicago, Chicago, Illinois 60637, USA}
\email{souganidis@math.uchicago.edu}
\thanks{P. Souganidis was partially supported by the National Science
Foundation. P. Cardaliaguet was partially supported by the ANR (Agence Nationale de la Recherche) through MICA project (ANR-06-BLAN-0082)}
\dedicatory{Version: \today}

\begin{document}

\begin{abstract}
We consider the so-called $G$-equation, a level set Hamilton-Jacobi equation, used as a sharp interface model for flame propagation,
perturbed by an oscillatory advection in a spatio-temporal periodic environment. Assuming that the advection
has suitably small spatial divergence, we prove that, as the size of the oscillations diminishes, the solutions
homogenize (average out) and converge to the solution of an effective anisotropic first-order (spatio-temporal
homogeneous) level set equation. Moreover we  obtain a rate of convergence
and show that, under certain conditions, the averaging enhances the velocity of the underlying front. We also prove that, at scale one,
the level sets of the solutions of the oscillatory problem converge,
at long times, to the Wulff shape associated with the effective Hamiltonian. Finally we also consider advection depending on position
at the integral scale.
\end{abstract}
%\tableofcontents

 \maketitle

\section{Introduction}
We study the limit, as $\ep\to0$, of the solution to the level-set equation
\be\label{eq:uep}
\left\{\begin{array}{l}
(i)\quad u_t^\ep=|Du^\ep|+\lg V(\frac{x}{\ep}, \frac{t}{\ep}),Du^\ep \rg  \quad {\rm in } \quad \R^N\times (0,T)\\
%\ee
%with initial condition
%\be\label{eq:ic}
(ii) \quad u^\epsilon = u_0 \quad \text { on } \quad \R^N\times\{0\}\ .
\end{array}\right.\ee

Equation \eqref{eq:uep}(i) is referred to as the $G$-equation, and is used as a model for flame propagation in turbulent fluids (\cite{Pet,Wlms85}). In that setting, the level sets of the function $u^\epsilon$ represent the evolving flame surface and $-V$ is the fluid velocity field.  At points where $u^\epsilon$ is differentiable and $| Du^\epsilon | \neq 0$, the level sets of $u^\epsilon$ move with normal velocity
$$\nu = 1 - \lg V(\frac{x}{\ep}, \frac{t}{\ep}), \hat n\rg \ ,$$
where $\hat n = -Du^\epsilon/|Du^\epsilon|$ is the exterior normal vector of the front. When $V \equiv 0$, level sets move with constant speed $s_L = 1$, which is called the laminar speed of flame propagation.

We assume that
the vector field $V\in {\mathcal C}^{0,1} (\R^{N+1}; \R^N)$ \quad is $\mathbb{Z}^{N+1}$-periodic in both $x$ and $t$, i.e., for all $(x,t) \in \R^{N+1}$,
$k \in \mathbb{Z}^N$ and $s \in \mathbb{Z}$,
\be\label{eq:periodic}
V(x+k,t + s) = V(x,t).
\ee

Our first result says that there exists a positively homogeneous of degree one, convex and continuous Hamiltonian $\bar H$ such that, as $\epsilon \to 0$, the $u^\ep$'s converge locally uniformly in $\R^N \times [0,\infty)$ to the solution $\bar u$ of the initial
value problem

\be\label{eq:hom}
\begin{cases}
\bar u_t = \bar H (D \bar u) \quad \text{ in } \quad  \R^N \times [0,\infty),\\
\noalign{\vskip6pt}
\bar u = u_0\ \text{ on }\ \R^N\times \{0\}\ .
\end{cases}
\ee

Although $V$ is bounded, we do not assume that $|V| < 1$, and, hence, the Hamiltonian $$H(x,t,p) = |p| + \langle V(x,t), p\rangle$$ is not coercive in $\abs{p}$ at every point $(x,t)$. This lack of coercivity is the main mathematical challenge in the analysis.  If either $|V|<1$ or the nonlinearity  $|Du|$ were replaced with $| Du |^2$, then $H$ would be coercive in $\abs{p}$ and the problem would be within the scope of the theory developed in \cite{LPV}. There are, however, relatively few homogenization results about noncoercive Hamiltonians \cite{AI01, ABMemoirs, GBa, Ca, IM}, and
none of them deals with the particular structure considered here.

The following simple example shows that, in the absence of coercivity, some additional assumption about the divergence of $V$ is necessary in order for the   $u^\epsilon$'s to have a local uniform limit. To this end, let $V = V(x)$ be a smooth $\mathbb{Z}^{N}$-periodic vector field such that, in the cube $Q_1=(-\frac12,\frac12)^N$,  $V(x) = -10 x$, if $\abs{x} < 1/6$, and $V(x) = 0$ if $\abs{x} \geq 1/3$.
It is known that the $u^\epsilon$'s  have the control representation
\be
u^\epsilon(x,t) = \sup (u_0(X_x(t))),
\ee
where the supremum is over all functions $X_x \in W^{1,\infty}([0,t];\R^N)$ such that $X_x(0) = x$ and $ X_x'(s) = \kappa(s) + V(X_x(s))$ with the controls $\kappa(\cdot)$ satisfying $|\kappa|\leq 1$. If $u_0(x) = \langle p, x \rangle$ with $\abs{p} > 0$, we see easily that $\lim_{\ep \to 0} \abs{u^\epsilon(0,t)} = 0$ for all $t > 0$. However, $\liminf_{\epsilon \to 0} u^\epsilon(x_\epsilon, t) > 0$ if $t > 0$ and $\{ x_\epsilon\}_\epsilon$ is any point satisfying $\abs{x_\epsilon} = \epsilon/2$.  Roughly speaking, the problem with such a vector field $V$ is that it traps the trajectories which start at the lattice points. %$x \in \mathbb{Z}^N$.
If the divergence of $V$ is sufficiently small, however, it is reasonable to expect that the controls are strong enough to overcome such traps.

We assume that $V$ has ``small divergence", in the sense that, for all $t\in\R$,
\be\label{divG}
\alpha(t) = \frac{1}{c_I} - \|{\rm div}_x V(\cdot,t) \|_{L^N(Q_1)}\ \geq 0 \quad \text{ and } \quad  \alpha^* = \int_0^1 \alpha(s) \,ds > 0,
\ee
%and $\alpha^* = \int_0^1 \alpha(s) \,ds > 0$.
where $c_I$ is the isoperimetric constant in the cube $Q_1$ (see, for instance, \cite{eg92}), i.e.,  the smallest constant such that, for all measurable subsets $E$ of $Q_1$,
$$
\left(|E|\wedge |Q_1\backslash E|\right)^{(N-1)/N} \leq c_I {\rm Per} (E, Q_1),
$$
and also the optimal constant for the Poincar\'e inequality
$$
\| f-\lg f\rg\|_{L^{1^*}(Q_1)} \leq c_I  \|Df\|_{L^1(Q_1)},
$$
for $f \in W^{1,1}(Q_1)$, $1^*=N/(N-1)$ and $\lg f\rg = \int_{Q_1}f(x)\,dx$.

To state the main results we introduce some additional notation. Throughout the paper we use $Q_1^+$ and $BUC(\bar U)$ to denote respectively the space-time cube $Q_1^+ = Q_1 \times [0,1] \subset \R^{N+1}$ and the space of bounded uniformly continuous functions on $\bar U$, and we write
$$
\lg V\rg =\int_{Q_1^+} V(x,t)dx\,dt \qquad {\rm and }\qquad \lg x {\rm div} V\rg = \int_{Q_1^+} x{\rm div} V(x,t)dx \,dt\;.
$$

We have:

\begin{Theorem}\label{theo:mainperiodic} Assume that $V\in {\mathcal C}^{0,1}(\R^{N+1};\R^N)$ satisfies \eqref{eq:periodic} and \eqref{divG}
and that $u_0\in {\mathcal C}^0(\R^N)$ is bounded. There exists a positively homogeneous of degree one,
Lipschitz continuous, convex Hamiltonian $\bar H:\R^N\to \R$ such that, if $u^\ep \in {\mathcal C}^0(\R^N\times[0,+\infty))$ and 
$\bar u \in {\mathcal C}^0(\R^N\times[0,+\infty))$
are the solutions  to the initial value problems (\ref{eq:uep}) and \eqref{eq:hom} respectively with initial datum $u_0$, then, as $\ep \to 0$, the  $u^\ep$'s converge locally uniformly in $\R^N\times[0,T]$ to $\bar u$.
Moreover,
%$\bar H(\gamma P) = \gamma \bar H(P)$ for all $\gamma > 0$
for all $P \in \R^N$,
\be\label{HPstrict}
\bar H(P)\ \geq |P| \int_0^1 \left(1-c_I\|{\rm div}V(\cdot,t)\|_{L^N(Q_I)}\right) \,dt+ \lg\ \lg V\rg +\lg x {\rm div} V\rg, P\rg.
%\ \rg \qquad \forall P\in \R^N\;.
\ee
Finally, the convex map $P\to \bar H(P) - \lg\ \lg V\rg +\lg x {\rm div} V\rg, P\ \rg$ is coercive.
\end{Theorem}

For Lipschitz continuous initial datum $u_0$, we can actually estimate the convergence rate as $\epsilon \to 0$. We have:

\begin{Theorem}\label{theo:error} Assume that  $u_0\in C^{0,1}(\R^N)$ and let $u^\epsilon, {\bar u} \in BUC(\R^N \times [0,T])$,
for all $T>0$, be respectively the solutions to (\ref{eq:uep}) and \eqref{eq:hom}. Then,
for all $T>0$, there exists a positive constant $C$ that depends only on $T$, $N$, $V$ and the Lipschitz constant of $u_0$, such that, for all
$(x,t)\in \R^N\times [0,T]$,
$$
|u(x,t)-u^\epsilon(x,t)|\leq C\epsilon^{1/3}.
$$
\end{Theorem}

In the case that, for all $x\in\R^N$ and $t\in\R$,
\be\label{eq:free}
{\rm div}_x V(x,t)=0,
\ee
we derive some additional properties of the function $\bar H$. To simplify the statement we
also assume, without any loss of generality
(see Lemma \ref{lossgen} below), that
\be\label{eq:mean}
%\lg V \rg =0
\int_{Q_1} V(x,t)dx=0\qquad \forall t\in \R.
 \ee

Then, according to Theorem \ref{theo:mainperiodic}, the averaged Hamiltonian $\bar H:\R^N\to\R$ satisfies, for all $P\in\R^N$,
$$
\bar H(P) \geq |P|.
$$

We establish here a necessary and sufficient condition to have the strict inequality
$$\bar H(P) > \abs{P},$$ in which case we have enhancement of the speed due to averaging.

Recall that, since $\bar H$ is homogeneous of degree one, the level sets of $\bar u$ move with speed $\bar H(\hat n)$ in the
direction of the normal vector $\hat n = -D \bar u/|D \bar u|$. Therefore, we refer to the situation $\bar H(P) > \abs{P}$ as ``enhancement", because it implies that such velocity fields lead to faster propagation of interfaces compared to the case $V \equiv 0$.

First, we state the result in the case where $V$ only depends on $x$. We have:

\begin{Theorem}\label{EnhenVx} Assume that $V\in C^{0,1}(\R^N)$ is $\Z^N$-periodic, ${\rm div}V=0$ and $\lg V\rg=0$, and
let $P\in\R^N\backslash \{0\}$. Then $\bar H(P)=|P|$ if and only if, for all $x\in \R^N$,  $\lg V(x), P\rg=0$.
In particular, if $N=2$, then $\bar H(P)=|P|$ if and only if the  stream function $E$ associated to $V$ is of the form $E=\tilde E(\lg P,\cdot\rg)$
for some $\tilde E:\R\to \R$, i.e., $V$ is a shear drift in the direction orthogonal to $P$.
\end{Theorem}

When $V$ is also time dependent, the characterization of equality $\bar H(P)=|P|$ is provided by

\begin{Theorem}\label{EnhenVxt} Assume that, for all $t\in\R$,  ${\rm div}_xV(\cdot,t)=0$ and $\int_{Q_1} V(x,t)dx=0$, and
fix $P\in\R^N\backslash \{0\}$. Then
$\bar H(P)=|P|$ if and only if
there exists $\hat z\in BV_{loc}(\R)$ such that $ \hat z'\geq -|P|$ in the sense
of distribution and the function $z(x,t)=\hat z(\frac{\lg P,x\rg}{|P|}+t)$
is $\Z^{N+1}-$periodic and satisfies, for all $t\in\R$, in the sense of distributions
\be\label{divzPV0}
{\rm div}\left( (z(\cdot,t)+\lg P, \cdot\rg) V(\cdot,t)\right)= 0 \quad{\rm in } \quad \R^{N}.
\ee
\end{Theorem}
	
We continue with  some observations about these results.  Theorem \ref{EnhenVx} yields that, if $N=2$, $\bar H(P)=|P|$ and $V$ is not constant,
then $P=(P_1,P_2)$ must be a rational direction, i.e., either $P_2= 0$ or $P_1/P_2\in \Q$, since
$V$ is $\Z^2-$periodic.  For Theorem \ref{EnhenVxt}, observe that, if $z$ is not constant, then
$P/|P|$ must be a rational vector.

Also (\ref{divzPV0}) is equivalent to saying (see Lemma \ref{constant} below) that, for any fixed $t>0$,
the map $x\to z(x,t)+\lg P,x\rg$ is constant along the flow of the differential equation
$X'(s)=V(X(s),t)$.

We remark that it is possible to construct nontrivial examples of time-dependent flows for which $\bar H(P) =|P|$.  Indeed when $N=2$ for any smooth,
$\mathbb{Z}^1$-periodic $(E_1,E_2)$ such that
$E_1(0)=0$, let $E(x_1,x_2,t)=E_1(x_1+t)E_2(x_2)$, $V= \nabla^\perp E$ and $P=(1,0)$. Then $\bar H(P)=|P|$ because the map $\hat z(s)=[s]-s$, where $[s]$
stands for the integer part of $s$, satisfies the condition of Theorem \ref{EnhenVxt}.  For more analysis and numerical computation of $\bar H$ for specific flow structures, we refer to \cite{It,EMS, MS,  NXG, Ob_01}.\\

The next result of the paper is about the long time behavior of the solution to \eqref{eq:uep} with $\ep=1$ and, in particular, the convergence, as
$t\to\infty$, of its
zero level set
to the Wulff-shape associated with the effective $\bar H$, which is given by
\be\label{eq:wulff}
{\mathcal W}=\{ y\in \R^N :\lg P,y\rg+\bar H(P)\geq 0\quad \text{ for all} \quad  P\in \R^N\}.
\ee

We consider the initial value problem
\be\label{eq:ep}
\begin{cases}
u_t=|Du|+\lg V\left(x, t \right),Du \rg  \quad {\rm in } \quad \R^N\times (0,\infty),\\
\noalign{\vskip6pt}
u = u_0\ \text{ on }\ \R^N\times \{0\}\ ,
\end{cases}
\ee
and set, for all $t\geq0$,

$$ K(t)=\{x\in\R^N : u(x,t)\geq0\}.$$

Recall that, in the language of front propagation (see, for example, \cite{bss93}), the family of closed sets $(K(t)_{t\geq0}$
is solution of the front propagation problem
$$  \nu=1-\lg V(x,t),\hat n \rg$$
starting from $K(0)=K_0$.

We have:

\begin{Theorem}\label{theo:wulff} Let $K_0$ be a non-empty  compact subset of $\R^N$.
%and $(K(t)))_{t\geq 0}$ be the solution to the front propagation problem
%$V_{x,t}=1-\lg V(x,t),\nu_x\rg$ starting from $K_0$.
There exist $C>0$ and $T>0$ such that, for all $t\geq T$,
\be\label{InclWulff1}
K(t)\ \subset \ (t+C) {\mathcal W}.
\ee
Moreover, there exists a constant $C_0>0$, independent of $K_0$, such that, if $K_0$ contains a cube of side length $C_0$,
then there exist $C>0$ and $T>0$ such that, for all $t\geq T$,
\be\label{InclWulff2}
(t-Ct^{2/3}) {\mathcal W} \ \subset K(t).
\ee
%for some constants $C>0$ and $T>0$.
\end{Theorem}

We note that we do not know whether the size condition on $K_0$ is actually necessary.

The final result of the paper is about homogenization when $V$ depends on $x$ also at the integral scale,
i.e., we are interested in the behavior as $\ep\to 0$ of the solutions to the initial value problem

\be\label{eq:xtyepsep}
\left\{\begin{array}{l}
u^\ep_t= |Du^\ep|+\lg V(x, \frac{x}{\ep},\frac{t}{\ep}), Du^\ep\rg \quad {\rm in }\quad \R^N\times (0,T)\\
u^\ep= u_0 \quad {\rm on } \quad \R^N\times\{0\}\ ,
\end{array}\right.
\ee

\noindent where $V:\R^{N}\times \R^N\times \R \to\R^N$ is smooth, bounded, $\Z^N$-periodic with respect to the last two
variables, i.e., for all $(x,y,s)\in\R^{N}\times \R^N\times\R$,
\be\label{eq:xper}
V(x,y+k,s+h)=V(x,y,s)\;,
\ee
is divergence free in the fast variable, i.e., for all $(x,y,s)\in\R^N\times \R^N\times\R$,
\be\label{eq:xdiv}
{\rm div}_yV(x,y,s)=0 \;,
\ee
and satisfies, for all $x\in\R^N$ the ``smallness'' condition
\be\label{BoundIntV}
\left|\int_0^1\int_{Q_1} V(x,y,s)dyds\right|<1\;.
\ee

The homogenized initial value problem is

\be\label{eq:HomPbVxxep}
\left\{\begin{array}{l}
\bar u_t= \bar H(x,D\bar u) \quad {\rm in } \quad \R^N\times (0,T),\\
\bar u= u_0 \quad {\rm on } \quad \R^N .
\end{array}\right.
\ee

We have:

\begin{Theorem}\label{HomVxxep} Assume \eqref{eq:xper}, \eqref{eq:xdiv} and \eqref{eq:HomPbVxxep}. There exists $\bar H \in {\mathcal C}^0(\R^N\times \R^N)$, which is positively homogeneous of degree one and convex
with respect to the second variable, such that,
for any initial condition $u_0\in BUC(\R^N)$, the solution $u^\ep$  to (\ref{eq:xtyepsep})
converges, as $\ep\to0$, locally uniformly in $\R^N\times[0,T]$ to the solution $\bar u$ of \eqref{eq:HomPbVxxep}.
Moreover $\bar H$ satisfies, for all $(x,P)\in \R^N\times \R^N$,
\be\label{CroissVxxep}
\bar H(x,P) \geq |P|+\lg \int_0^1\int_{Q_1} V(x,y,s)dyds , P\rg \;.
\ee
\end{Theorem}

The paper is organized as follows. Theorem \ref{theo:mainperiodic} and Theorem \ref{theo:error} are proved in Section \ref{sec:homogen}. In Section \ref{sec:enh} we prove Theorem \ref{EnhenVx} and Theorem \ref{EnhenVxt}, while Theorem\ref{theo:wulff} is proved in Section \ref{sec:Wulff}. In Section \ref{sec:xtys}, we prove an extension of Theorem \ref{theo:mainperiodic} to the case where $V = V(x,x/\ep,t/\ep)$ has large-scale spatial variation.
%In Section \ref{sec:future} we discuss some open problems.
The Appendix contains a proof of Lemma \ref{EstiPeri}, which plays an important role in the proof of Theorem \ref{theo:mainperiodic}.

\vspace{0.1in}

About the time this paper was completed, we learned about a similar but less general homogenization result obtained by different methods in \cite{xy}. In
particular it is proved in \cite{xy} that homogenization takes place for time independent advection satisfying $V = V_1 + V_2$ with ${\rm div } V_1 = 0$ and $| V_2| < 1$.

Finally we remark that throughout the paper we will need some basic results from the theory of viscosity solutions, like comparison principles, representation formulae, etc.. All such results can be found, for instance, in \cite{b}.

%%%%%%%%%%%%%%%%%%%%%%%%%%%%%%%%%%%%%%%%%%%%%%%%%%%%%%%
%%%%%%%%%%%%%%%%%%%%%%%%%%%%%%%%%%%%%%%%%%%%%%%%%%%%%%%
\section{Homogenization} \label{sec:homogen}

We begin with some preliminary discussion and results to set the necessary background for the proofs of Theorem \ref{theo:mainperiodic} and Theorem \ref{theo:error}. First, we recall that for any $\lambda > 0$ and any $P \in \R^N$,  the ``penalized'' cell problem
\be\label{eq:vlambda}
v_{\lambda,t} + \lambda v_\lambda =|Dv_\lambda +P|+\lg V, Dv_\lambda+P\rg  \quad \text{ in } \quad \R^{N+1},
\ee
has a unique $\mathbb{Z}^{N+1}$-periodic solution $v_\lambda \in BUC(\R^{N+1})$, which is actually H\"{o}lder continuous and satisfies, for all $(x,t)\in\R^{N+1}$, the bound
\be\label{eq:lbound}
- \lambda^{-1}\abs{P}(1 + \norm{V}_\infty) \leq v_\lambda  \leq  \lambda^{-1} \abs{P}(1 + \norm{V}_\infty).
\ee

We also recall that, in the periodic setting, homogenization is equivalent to proving that the
$(\lambda v_\lambda)$'s converge uniformly in $\R^{N+1}$, as $\lambda \to0$,  to some constant $\bar c(P)$ and that $\bar c \in {\mathcal C}^0(\R^N)$.
In this case, $\bar H(P)=\bar c(P)$. \\

In view \eqref{eq:lbound}, to prove the convergence of the $\lambda v_\lambda$'s, we need to control  their oscillations. We have:

\begin{Lemma}\label{osc} For all $P \in \R^N$ and $\lambda \in (0,1]$,
\be\label{eq:osc}
{\rm osc} (\lambda v_\lambda) \leq C|P| \lambda,
\ee
where
\[
C = 4(1+\|V\|_{\infty}) ( { 2^{1/N} \alpha^* }^{-1}N  + 3).
\]
\end{Lemma}

Before we present the proof of Lemma \ref{osc}, which is one of the most important parts of the paper, we
point out its main consequence in the next

\begin{Corollary}\label{cor:barH} Let $C$ be the constant given by Lemma \ref{osc}. There exists some $\bar H(P)\in  \R$ such that
$$
\|\lambda v_\lambda -\bar H(P)\|_\infty\leq C|P| \lambda\;.
$$
\end{Corollary}

\noindent {\bf Proof:}
%of Corollary \ref{cor:barH}: }
The maps $\lambda \to \lambda \min v_\lambda$ and $\lambda \to \lambda \max v_\lambda$ are respectively nonincreasing
and nondecreasing. For the sake of completeness we present a formal proof, which can be easily made rigorous using viscosity solutions
arguments. Since the two claims are proved similarly, we present details only for the former.

To this end, for $0<\lambda <\mu$, let
$(x,t)$ be a maximum of $v_\mu-v_\lambda$. Then,
at least formally, at $(x,t)$, we have $D_{x,t}v_\mu=D_{x,t} v_\lambda$, and
$$
v_{{\mu},t} +\mu v_\mu \leq |Dv_\mu +P|+\lg V, Dv_\mu+P\rg
$$
and
$$
v_{{\lambda},t} +\lambda v_\lambda \geq |Dv_\lambda+P|+ \lg V,Dv_\lambda+P\rg\;.
$$

It follows that $$\mu v_\mu(x,t) \leq \lambda v_\lambda(x,t)$$.

Let $(y,s)$ be a minimum point of $\lambda v_\lambda$. Then
$$
\begin{array}{rl}
\mu v_\mu(y,s) \; \leq &  \mu (v_\mu(x,t)-v_\lambda(x,t)+v_\lambda(y,s) ) \\
 \leq  & \lambda v_\lambda(x,t)-\mu v_\lambda(x,t)+\mu v_\lambda(y,s) \\
  \leq  & \lambda v_\lambda(x,t)+ \lambda ( v_\lambda(y,s)- v_\lambda(x,t)) \;  \leq \; \lambda  v_\lambda(y,s),
  \end{array}
$$
and, hence, $$\min \mu v_\mu \leq \min \lambda v_\lambda\ .$$

The above remark combined with Lemma \ref{osc} implies that the $\lambda v_\lambda$'s
converge uniformly to some constant $\bar H(P)$ and that
$$
\lambda \min v_\lambda \leq \bar H(P) \leq \lambda \max v_\lambda\;.
$$
%Using the bound on the oscillation gives the result.
\QED

We continue with the

\noindent {\bf Proof of Lemma \ref{osc}: } Without any loss of generality, we may assume that $V$ is smooth. Indeed, if the result
holds for any smooth $V$, then it also holds by approximation for any $V\in C^{0,1}$.

Recalling \eqref{eq:lbound}, $w_\lambda=\lambda v_\lambda$ satisfies, in the viscosity sense,
\be\label{ineqwlambda}
w_{\lambda,t} -|Dw_\lambda(x,t)|-\lg V(x,t),Dw_\lambda(x,t)\rg \geq -C_0\lambda \qquad {\rm in }\;\R^{N+1} \;,
\ee
where $$C_0=2(1+\|V\|_\infty)|P|\ .$$

It follows that $(x,t)\to w_\lambda (x,t) +C_0\lambda t$ is a viscosity super-solution of the level-set initial value problem
%\be\label{evolz}
%\left\{\begin{aligned}
%& (i) & \partial_t z(x,t)=|Dz(x,t)|+\lg V(x,t) , z(x,t))\rg  \qquad & {\rm in }\;\R^N \times (0,T)\\
%& (ii) & z(x,0)=w_\lambda (x,0) \qquad \qquad & {\rm in }\; \R^N
%\end{aligned}\right.
%\ee
\be\label{evolz}
\begin{cases}
(i) \quad z_t =|Dz|+\lg V , Dz\rg \quad {\rm in } \quad \R^N\times (0,\infty),\\
\noalign{\vskip6pt}
(ii) \quad z = w_{\lambda} \quad \text{ on } \quad  \R^N\times\{0\}\ .
\end{cases}
\ee

The standard comparison of viscosity solutions then implies that, for all $(x,t)\in\R^{N+1}$,
\be\label{comparaison}
w_\lambda (x,t)+C_0\lambda t \geq z(x,t)\;.
\ee
%where $z$ is the solution to (\ref{evolz}).

The next step is to understand the evolution of the perimeter of the level-sets of $z$. For this we need the following
result, which is proved in the Appendix.

We have:

\begin{Lemma}\label{EstiPeri} Assume that $V \in C^{1,1}(\R^{N+1})$ and let $z \in BUC(\R^{N+1})$ be a solution of
%$V:\R^{N+1}\to \R^N$ is of class ${\mathcal C}^{1,1}$ and $z$ is a continuous solution of
(\ref{evolz})(i).  Then, for any level
$\theta\in (\ \inf z(\cdot,0), \sup z(\cdot,0)\ )$ such that 
$$
\{z(\cdot,0)=\theta\}=\partial \{z(\cdot,0)>\theta\}=\partial \{z(\cdot,0)<\theta\}
\quad {\rm and }\quad |\{z(\cdot,0)=\theta\}|=0\;,
$$
and for all $t>0$, we have 
$$
\{z(\cdot,t)=\theta\}=\partial \{z(\cdot,t)>\theta\}=\partial \{z(\cdot,t)<\theta\}
$$
as long as $\{z(\cdot,t)<\theta\}\neq \emptyset$. Moreover the sets $\{z(\cdot,t)>\theta\}$ and $\{z(\cdot,t)<\theta\}$ have locally of finite perimeter, and $|\{z(\cdot,t)=\theta\}|=0$.
%the set $\{z(\cdot,t)>\theta\}$ satisfies:
%\begin{itemize}
%\item $\{z(\cdot,t)=\theta\}=\partial \{z(\cdot,t)>\theta\}=\partial \{z(\cdot,t)<\theta\}$,
%\item the sets $\{z(\cdot,t)>\theta\}$ and $\{z(\cdot,t)<\theta\}$ are locally of finite perimeter,
%\item for any continuous map $\varphi:\R^N\to \R$ with compact support, the map
%$\displaystyle{ I(t):= \int_{ \{z(\cdot,t)>\theta\}}\varphi(x)dx  }$ is absolutely continuous and satisfies
%$$
%\frac{d}{dt} I(t) = \int_{\{z(\cdot,t)=\theta\}} \varphi(x)(1-\lg V(x,t), \nu(x,t)\rg)d{\mathcal H}^{N-1}(x)
%$$
Finally, for any compactly supported $\varphi \in {\mathcal C}^0(\R^N)$, the maps
$t\to I(t)= \int_{ \{z(\cdot,t)>\theta\}}\varphi(x)dx  $ and $t\to J(t)=\int_{ \{z(\cdot,t)<\theta\}}\varphi(x)dx $ are
absolutely continuous and satisfy, for almost all $t>0$,
$$
\frac{d}{dt} I(t) = \int_{\{z(\cdot,t)=\theta\}} \varphi(x)(1-\lg V(x,t), \nu(x,t)\rg)d{\mathcal H}^{N-1}(x),
$$
and
$$
\frac{d}{dt} J(t) = -\int_{\{z(\cdot,t)=\theta\}} \varphi(x)(1+\lg V(x,t), \nu(x,t)\rg)d{\mathcal H}^{N-1}(x),
$$
where $\nu(x,t)$ denotes in the former identity the measure theoretic outward
unit normal to $\{z(\cdot,t)>\theta\}$ at $x\in \partial \{z(\cdot,t)>\theta\}$, while in the latter is
the measure theoretic outward
unit normal to $\{z(\cdot,t)<\theta\}$ at $x\in \partial \{z(\cdot,t)<\theta\}$.
%\end{itemize}
\end{Lemma}

%\begin{Remark}{\rm Since the set $\{z(\cdot,t)=\theta\}$ has a zero Lebesgue measure, we also have, for any
%continuous map $\varphi:\R^N\to \R$ with compact support, that the map
%$I(t):=\displaystyle{ \int_{ \{z(\cdot,t)<\theta\}}\varphi(x)dx  }$ is absolutely continuous and satisfies
%$$
%\frac{d}{dt} I(t) = -\int_{\{z(\cdot,t)=\theta\}} \varphi(x)(1+\lg V(x,t), \nu(x,t)\rg)d{\mathcal H}^{N-1}(x)
%$$
%for almost all $t>0$, where $\nu(x,t)$ is the measure theoretic outward
%unit normal to $\{z(\cdot,t)<\theta\}$ at $x\in \partial \{z(\cdot,t)<\theta\}$.
%}\end{Remark}

Continuing with the ongoing proof, suppose that there exists $\theta\in \R$ with
$$
\inf w_\lambda(\cdot, 0)< \theta<  \sup w_\lambda(\cdot,0)\;,
$$
such that
$$
\{w_\lambda(\cdot,0)=\theta\}=\partial \{w_\lambda(\cdot,0)>\theta\}=\partial \{w_\lambda(\cdot,0)<\theta\}
\quad {\rm and }\quad |\{w_\lambda(\cdot,0)=\theta\}|=0\;,
$$
and such that 
$$
\abs{ \{ x \in Q_1 : w_\lambda(x,0) < \theta \} } < 1/2\ .
$$
For all $t\geq0$, set
$$
\rho(t)=\left|\{ x \in Q_1 \;:\; z(x,t)<\theta\}\right| \;.
$$

Let $[0,T)$ be the maximal interval on which $\rho(t)<1/2$ for any $t\in [0,T)$. Note that  $T>0$ because
$\rho(0)<1/2$.  We claim that, for all $0\leq t_1\leq t_2<T$,
\be\label{EstiDecroissance}
\rho(t_2)-\rho(t_1) \leq -\int_{t_1}^{t_2}\alpha(s) \rho(s)^{(N-1)/N}ds \;.
\ee

Indeed
%\noindent {\it Proof of (\ref{EstiDecroissance}):}
fix a positive integer $R$ and let $Q_R=(-\frac{R}{2}, \frac{R}{2})^N$.
%\in \N^*$
The space periodicity of $z$ gives
$$
\rho(t)=R^{-N} \left|\{ x \in Q_R \;:\; z(x,t)<\theta\}\right| \;.
$$
%where $Q_R=(-\frac{R}{2}, \frac{R}{2})^N$.

For $h>0$ small, let $\chi_h\in C^{0,1} (\R^N;[0,1])$ be  such that
$$
\chi_h=1 \; {\rm in }\; Q_{R}\qquad {\rm and }\qquad \chi_h=0 \; {\rm in }\; \R^N\backslash Q_{R+h}\;,
$$
and, for any $t\in [0,T]$, set
$$
\rho_{R,h}(t)=R^{-N} \int_{\{z(\cdot,t)<\theta\}} \chi_h(x)dx\;,
$$
and note that,
$$
\lim_{h\to0} \rho_{R,h}(t)= \rho(t)\;.
$$

It follows from Lemma \ref{EstiPeri} that, for almost all $t\in (0,T)$,
$$
\frac{d}{dt} \rho_{R,h}(t)=-R^{-N} \int_{\{z(\cdot,t)=\theta\}} \chi_h(x)(1+\lg V(x,t), \nu(x,t)\rg)d{\mathcal H}^{N-1}(x)\;.
$$

Moreover
$$
\int_{\{z(\cdot,t)=\theta\}} \chi_h(x) d{\mathcal H}^{N-1}(x)\geq
{\mathcal H}^{N-1}(\{z(\cdot,t)=\theta\}\cap Q_R),
$$
and, in view of the spatial periodicity of $z$,
$$
{\mathcal H}^{N-1}(\{z(\cdot,t)=\theta\}\cap Q_R)
\geq R^N {\mathcal H}^{N-1}(\{z(\cdot,t)=\theta\}\cap Q_1)\;.
$$

The isoperimetric inequality in the box $Q_1$ and the fact that
$\left|\{z(\cdot,t)<\theta\}\cap Q_1\right|<1/2$ give
$$
{\mathcal H}^{N-1}(\{z(\cdot,t)=\theta\}\cap Q_1)
\geq \frac{1}{c_I} \left|\{z(\cdot,t)<\theta\}\cap Q_1\right|^{(N-1)/N}.
$$

Using once more the space periodicity of $z$ we get
$$
\begin{aligned}
& \left|\{z(\cdot,t)<\theta\}\cap Q_1\right|^{(N-1)/N}
&\; =\; & (R+1)^{-(N-1)}\left|\{z(\cdot,t)<\theta\}\cap Q_{R+1}\right|^{(N-1)/N} \\
&&\; \geq\;
%& \left(\frac{R}{R+1}\right)^{N-1}
& (R(R+1)^{-1})^{N-1} \left(\rho_{R,h}(t)\right)^{(N-1)/N}.
\end{aligned}$$

Combining all the above we obtain
$$
R^{-N}\int_{\{z(\cdot,t)=\theta\}} \chi_h(x) d{\mathcal H}^{N-1}(x)\geq
\frac{1}{c_I}
%\left(\frac{R}{R+1}\right)^{N-1}
(R(R+1)^{-1})^{N-1}\left(\rho_{R,h}(t)\right)^{(N-1)/N}\;.
$$

Next we estimate the integral
$$\displaystyle{  \int_{\{z(\cdot,t)=\theta\}} \chi_h(x)\lg V(x,t), \nu(x,t)\rg d{\mathcal H}^{N-1}(x) }.$$

For some constant $k$ depending only on $N$ we  have:
$$
\begin{aligned}
& -\int_{\{z(\cdot,t)=\theta\}} \chi_h(x)\lg V(x,t), \nu(x,t)\rg d{\mathcal H}^{N-1}(x) \\
&  = - \int_{\{z(\cdot,t)<\theta\} } {\rm div}(\chi_h V)(x,t)\ dx\\
& \leq \;    -\int_{\{z(\cdot,t)<\theta\} } \chi_h(x) {\rm div}V(x,t)dx+ \|D\chi_h\|_\infty \|V\|_\infty \left|Q_{R+h}\backslash Q_{R}\right|\\
& \leq \;   \ (\int_{Q_{R+1} } \left|{\rm div}V(x,t)\right|^N\ dx)^{1/N}
(\int_{\{z(\cdot,t)<\theta\} } {\chi_h(x)}^{N/(N-1)} dx)^{(N-1)/N}\\
& \quad \qquad \qquad \qquad +  k R^{N-1}\|V\|_\infty \\
& \leq \;
%&\qquad\leq \;
(R+1)R^{N-1} \|{\rm div}V(\cdot,t)\|_{L^N(Q_1)}\left(\rho_{R,h}(t)\right)^{(N-1)/N}+ 
%\\
%& \qquad 
k R^{N-1}\|V\|_\infty .
\end{aligned}
$$

Hence, for almost all $t\in (0,T)$, we get
$$
\frac{d}{dt}\rho_{R,h}(t)\leq - (\frac{1}{c_I}\left(\frac{R}{R+1}\right)^{N-1}- \frac{R+1}{R}\|{\rm div}V(\cdot,t)\|_{L^N(Q_1)})
  (\rho_{R,h}(t))^{(N-1)/N}+ \frac{k}{R}\|V\|_\infty\;.
$$

Integrating first over $[t_1,t_2]$ and then letting $h\to0$ and $R\to +\infty$ gives (\ref{EstiDecroissance}). \\

Since $\alpha(t) \geq 0$, $\rho(t)$ is non increasing on $[0,T)$. Hence $T=+\infty$. Integrating (\ref{EstiDecroissance}) over $(0,t)$ we obtain,
for every $t\geq0$, that
$$
\rho(t)\leq \left[\rho^{1/N}(0)-\frac{1}{N} \int_0^t \alpha(s) \,ds \right]^N_+ \;,
$$
where $[s]_+=\max\{s,0\}$.

From the assumption $\rho(0)< 1/2$, it follows that
\[
t^* =  1 + \frac{N}{2^{1/N}\alpha^*} \geq \inf \left\{ t \;:\;  \int_0^t \alpha(s) \,ds \geq \frac{N}{2^{1/N}} \right\},
\]
and, hence, $\rho=0$ in $[t^*,\infty)$. The
continuity and the spatial periodicity of $z$ then yield that
$$z\geq \theta \quad \text{ in } \quad \R^N \times [t^*,\infty).$$

Let $k$ be an integer in the interval $[t^*, t^* + 1]$. The space-time periodicity of $w_\lambda$ and (\ref{comparaison}) give
\br
\inf_{t \in [0,1]} \inf_{x \in \R^N} w_\lambda(x,t) & = & \inf_{t \in [k,k + 1]} \inf_{x \in \R^N} w_\lambda(x,t) \nonumber \\
& \geq & \inf_{t \in [k,k + 1]} \inf_{x \in \R^N} z(x,t) - C_0 \lambda (t^* + 2) \nonumber \\
& \geq & \theta - C_0 \lambda (t^* + 2) .
\er

It follows that, if $\theta \in \R$ is such that
\be
\abs{ \{ x \in Q_1 : w_\lambda(x,0) < \theta\} } < 1/2,
\ee
then
\be
\inf_{t \in [0,1]} \inf_{x \in \R^N} w_\lambda(x,t) \geq \theta - C \lambda \label{EstiInf}
\ee
where $$C =  C_0(t^* + 2).$$

%%%%%%%%%%
\commentout{
Using (\ref{comparaison}) we obtain that
$$
w_\lambda(x)+C_0\lambda t^* \geq \theta\qquad \forall x\in \R^N.
$$
Since $C\geq C_0 t^*$,  we have proved that, for any $\theta\in (\min w_\lambda, \max w_\lambda)$ such that $|\{w_\lambda <\theta\}\cap Q_1|< \frac12$, we have
\be\label{EstiInf}
w_\lambda\geq \theta-C\lambda\;.
\ee
}%%%%%%%%%%%%%%%%%%%end commentout

Now we derive an upper bound. To this end suppose that $\theta \in \R$ 
with
$$
\inf w_\lambda(\cdot, 0)< \theta<  \sup w_\lambda(\cdot,0)\;,
$$
such that
$$
\{w_\lambda(\cdot,0)=\theta\}=\partial \{w_\lambda(\cdot,0)>\theta\}=\partial \{w_\lambda(\cdot,0)<\theta\}
\quad {\rm and }\quad |\{w_\lambda(\cdot,0)=\theta\}|=0\;,
$$
and such that 
\be
\abs{ \{ x \in Q_1 : w_\lambda(x,0) > \theta \} } < 1/2 .
\ee
The claim is that
\be\label{EstiSup}
\max_{x,t} w_\lambda(x,t) \leq \theta+C\lambda\;,
\ee
where $C=2|P|(1+\|V\|_\infty)(2+N2^{-1/N}(\alpha^*)^{-1})$. 
Indeed, arguing by contradiction, we assume that $$\max_{x,t} w_\lambda(x,t)> \theta+C\lambda\;.$$
Then, by continuity and periodicity of $w_\lambda$, there is some $\tau\in [0,1]$ such that 
$$
\abs{ \{x\in Q_1\ : \ w_\lambda(x,\tau)>\theta+C\lambda\} } >0\ .
$$ 
Let $z$ satisfy (\ref{evolz}-(i)) with initial condition $z=w_\lambda$ on $\R^N\times \{\tau \}$. As before, we have, for all $(x,t)\in\R^{N}\times[\tau, +\infty)$,
$$w_\lambda(x,t) \geq z(x,t) - C_0 \lambda (t-\tau).$$
Set
\be
\rho(t) = \abs{ \{ x \in Q_1 \; : \; z(x,t) > \theta + C \lambda \} }\qquad t\geq \tau\ ,
\ee
and observe that $\rho$ is continuous with $\rho(\tau)>0$. Then, arguing as before, we find that, for all $\tau\leq t_1\leq t_2$,
\be\label{eq:rho}
\rho(t_2)-\rho(t_1)\geq  \int_{t_1}^{t_2} \alpha(s) \left(\min\{\rho(s), 1-\rho(s)\}\right)^{(N-1)/N}ds \;.
\ee
Since $\alpha(t) \geq 0$ for all $t\geq0$, %\eqref{eq:rho} implies
it follows that $\rho$ is nondecreasing on $[\tau,+\infty)$. 
We claim that there is some $T\leq t^* = \tau+1 + N/ (2^{1/N} \alpha^*)$ such that 
$\rho\geq 1/2$ for $t\geq T$. Indeed, otherwise, $\rho < 1/2$ on $[\tau,t^*]$ and integrating \eqref{eq:rho} over $[\tau,t^*]$ gives
$$
\rho(t^*) \geq (\frac{1}{N} \int_\tau^{t^*} \alpha(s)\ ds)^N \ >\ \frac12\;,
$$
which contradicts our assumption. Let now $k$ be an integer in $[T,T + 1]$. The time-periodicity of $w_\lambda$ yields
$$
%\br
\begin{array}{rl}
 1/2 \;\leq   & \ds{ \rho(k) = \abs{ \{ x \in Q_1 \;:\; z(x,k ) > \theta + C \lambda\} }  }\\
\leq & \ds{ \abs{ \{ x \in Q_1 \;:\; w_\lambda(x,k ) + C_0 \lambda (k-\tau) > \theta + C \lambda \} }  }\\
= & \ds{ \abs{ \{ x \in Q_1 \;:\; w_\lambda(x,0) + C_0 \lambda (k-\tau) > \theta + C \lambda \} }.  }
\end{array}
%\er
$$
It follows from $ k \leq (t^* + 1)$ that
$$C_0 (k-\tau) \leq C \quad \text{ and }
\abs{ \{ x \in Q_1 \;:\; w_\lambda(x,0)  > \theta  \} } \geq 1/2,
$$
which contradicts the definition of $\theta$, and, hence, (\ref{EstiSup}) holds.
%%\be
%%\max_x w_\lambda(x,0) \leq \theta + C \lambda.
%%\ee
%
%Since $(x,t)\to w_\lambda(x,t) - C_0 \lambda t$ is a subsolution to (\ref{evolz}), the time-periodicity of $w_\lambda$ and the comparison principle imply that
%\br \label{tupper2}
%\max_{x \in \R^N} \max_{t \in [0,1]} w_\lambda(x,t) & = & \max_{x \in \R^N} \max_{t \in [1,2]} w_\lambda(x,t)  \nonumber \\
%& \leq & \max_{x \in \R^N}  \max_{t \in [1,2]} \left( w_\lambda(x,0) + C_0 \lambda t \right) \nonumber \\
%& \leq & \theta + C \lambda + 2 C_0 \lambda  \leq \theta +  2 \lambda C.
%\er

Finally set
\be
\bar \theta = \sup \left \{ \theta \in \R \;:\; \;\abs{ \{x \in Q_1\,:\, w_\lambda(x,0) < \theta \} } < \frac{1}{2}\right \} \;.
\ee
In view of the above, using (\ref{EstiInf}) and (\ref{EstiSup}), we get
\be
\min_{t \in [0,1]} \min_{x \in \R^N} w_\lambda(x,t) \geq \bar \theta - C \lambda
\ee
and
\be
\max_{t \in [0,1]} \max_{x \in \R^N} w_\lambda(x,t) \leq \bar \theta +  C \lambda\;,
\ee
where $C=2|P|(1+\|V\|_\infty)(3+N2^{-1/N}/\alpha^*)$. 
It follows that $\text{osc}(w_\lambda) \leq 2 C \lambda$, and, therefore, (\ref{eq:osc}).
\QED

We proceed with the

\noindent {\bf Proof of Theorem \ref{theo:mainperiodic} : } Let $\bar H(P)$ be defined by Corollary \ref{cor:barH}.
The fact that the map  $P\to \bar H(P)$ is positively homogeneous, convex and Lipschitz continuous follows easily from the properties of
(\ref{eq:vlambda}) and the comparison principle of viscosity solutions.\\

To prove (\ref{HPstrict}), first we perturb (\ref{eq:vlambda}) by a vanishing viscosity, i.e., for $\eta>0$ we consider
\be\label{eq:vlambdaeta}
v_{\lambda,t}^\eta + \lambda v_\lambda^\eta -\eta \Delta v_\lambda^\eta =|Dv_\lambda^\eta+P|+\lg V,Dv_\lambda^\eta+P\rg  \qquad {\rm in }\;\R^{N+1},
\ee
which has a unique $\Z^{N+1}-$periodic solution $v_\lambda^\eta \in BUC(\R^{N+1})$ which is at least in ${\mathcal C}^1(\R^{N+1})$ and
converges uniformly, as $\eta\to0$, to $v_\lambda$.

Integrating (\ref{eq:vlambdaeta}) over $Q_1^+$ and using the periodicity, we find
$$
\int_{Q_1^+} \lambda v_\lambda^\eta dx\,dt = \int_{Q_1^+} |Dv_\lambda^\eta +P|dx\,dt+ \int_{Q_1^+} \lg V,Dv_\lambda^\eta \rg dx\,dt +\lg \lg V\rg ,P\rg\;.
$$

Set
$$\xi(x)=\lg P, x\rg ,  \quad \lg v_\lambda^\eta(t)\rg= \int_{Q_1} v_\lambda^\eta (x,t)dx \quad \text{ and } \quad
\lg \xi\rg = \int_{Q_1} \xi(x)dx.$$
%and
%$$\quad \lg v_\lambda^\eta\rg_t= \int_{Q_1} v_\lambda^\eta (x,t)dx \quad \text{ and } \quad \lg xdivv_\lambda^\eta \lg_t =
%\int_{Q_1} x div v_\lambda^\eta (x,t)dx.
%$$

Since both $V$ and $v_\lambda$ are $\Z^{N+1}$-periodic, for each $t \in [0,1]$, we have
 $$
 \begin{aligned}
  &\int_{Q_1} \lg V(x,t),Dv_\lambda^\eta(x,t)\rg dx\\
  & \quad  =\; - \int_{Q_1} (v_\lambda^\eta(x,t)-\lg v_\lambda^\eta(t)\rg)\, {\rm div}V(x,t) \ dx\\
  & \quad  =\; - \int_{Q_1} (v_\lambda^\eta(x,t)-\lg v_\lambda^\eta(t)\rg +\xi(x)-\lg \xi\rg) {\rm div}V(x,t) \ dx+\\
  & \qquad\qquad\qquad\qquad \int_{Q_1} (\xi(x)-\lg \xi\rg){\rm div} V(x,t) dx.
  \end{aligned}
   $$

Applying, for each $t\in[0,1]$, H\"{o}lder's and Poincar\'e's inequalities yields
   $$
    \begin{aligned}
   & \int_{Q_1} (v_\lambda^\eta(x,t)-\lg v_\lambda^\eta\rg_t +\xi(x)-\lg \xi\rg) {\rm div}V(x,t) \ dx \\
    & \quad  \leq \; ( \int_{Q_1} |v_\lambda^\eta(x,t)-\lg v_\lambda^\eta\rg_t +\xi(x)-\lg \xi\rg|^{N/(N-1)}dx)^{(N-1)/N}(\int_{Q_1} |{\rm div}V(x,t)|^Ndx\ )^{1/N}\\
   &\qquad  \leq \;   c_I \  (\int_{Q_1}|Dv_\lambda^\eta(x,t)  +P|dx) \ \|{\rm div}V(\cdot, t)\|_{L^N(Q_I)}  \;,
  \end{aligned}
   $$
and, since
%$\displaystyle{ \int_{Q_1} (\xi(x)-\lg \xi\rg){\rm div} V(x,t) dx= \lg \ \lg x {\rm div }V\rg_t \ , \ P\rg }$,
$$
\int_{Q_1^+} (\xi(x)-\lg \xi\rg){\rm div} V(x,t) dxdt= \lg \ \lg x {\rm div }V\rg \ , \ P\rg ,
$$

\br\label{HPstrict1}
\int_{Q_1^+} \lambda v_\lambda^\eta (x,t)dx\,dt & \geq &\int_0^1 (1-c_I \|{\rm div}V(\cdot,t)\|_{L^N(Q_I)})
\int_{Q_1} |Dv_\lambda^\eta (x,t)+P|dx \,dt \nonumber \\
& & + \lg \lg V\rg + \lg x {\rm div }V\rg,P\rg\;.
\er

Finally, in view of (\ref{divG}), we have, for all $t\geq0$, $$1-c_I \|{\rm div}V(\cdot,t)\|_{L^N(Q_I)}\geq 0,$$ while the
periodicity of $v_\lambda$ implies, again for all $t\geq0$, that
$$
|P|= \left| \int_{Q_1} (Dv_\lambda^\eta (x,t)+P)dx\right| \leq \int_{Q_1} |Dv_\lambda^\eta (x,t)+P|dx\;.
$$
Letting $\eta\to0$ and then $\lambda\to0$ in (\ref{HPstrict1}) we obtain  (\ref{HPstrict}).
\QED

Error estimates for the periodic homogenization of coercive Hamilton-Jacobi equations were obtained earlier in  \cite{CDI2001}.
Although the proof of Theorem \ref{theo:error} is almost the same as the one of the analogous result in \cite{CDI2001},
we present it here for the sake of completeness. To simplify the presentation below we denote by $C$ constants that may change from line to line
but depend only on $u_0$, $\|Du_0\|_\infty$ and $V$.

We have:

\noindent {\bf Proof of Theorem \ref{theo:error} : } For all $(x,t,P)\in  \R^N\times \R\times \R^N$, set
$$
H(x,t,P)= |P|+ \lg V(x,t), P\rg \;.
$$

To avoid any technical difficulties due to the unboundedness of the domain, first we assume that $u_0$ is $(M\epsilon)\Z^N-$periodic for
some positive integer $M$, which implies that $u^\epsilon$ and $\bar u$ are
also $(M\epsilon)\Z^N-$periodic, and we obtain the estimate with constant independent of $M$. Then we use
the finite speed of propagation property
of the averaged initial value problem, to remove the restriction on $u_0$.

Let $v_\lambda=v_\lambda(\cdot,\cdot;P)\in BUC(\R^{N+1})$ be the $\Z^{N+1}$-periodic solution to \eqref{eq:vlambda}
%$$
%v_{\lambda,t} + \lambda v_\lambda=\left|Dv_\lambda+P\right|+\lg V, Dv_\lambda +P\rg\qquad {\rm in }\; \R^{N+1}\;,
%$$
and recall that the map
\be\label{LipVlambda}
\mbox{\rm $P\to \lambda v_\lambda(\cdot,\cdot;P)$ is $(1+\|V\|_\infty)$-Lipschitz continuous; }
\ee
to simplify statements heretofore we say the $f$ is $L$-Lipschitz continuous if it is Lipschitz continuous with
constant at most $L$.

Fix $T>0$ and consider  $\Phi:\R^N\times [0,T]\times \R^N\times [0,T]\to\R$ given by
$$
\Phi(x,t,y,s)=u^\ep(x,t)-u(y,s)-\ep v_\lambda\left(\frac{x}{\ep},\frac{t}{\ep},\frac{x-y}{\ep^\beta}\right)
-\frac{|x-y|^2}{2\ep^\beta} -\frac{(t-s)^2}{2\ep}-\delta s ,
$$
where $\beta\in (0,1)$, $\lambda\in(0,1)$ and $\delta>0$ are to be chosen later.

Since $\Phi$ is periodic in the space variables,
it has a maximum at some point $(\hat x,\hat t, \hat y,\hat s)$. The main part of the proof consists in showing that
either $\hat t=0$ or $\hat s=0$ for a suitable choice of $\lambda$ and $\delta$. \\

We argue by contradiction and assume that $\hat t>0$ and $\hat s>0$. Since $\bar u$ is Lipschitz continuous
and (\ref{LipVlambda}) holds, by standard arguments from the theory of viscosity solutions, we have
\be\label{Sti0}
\frac{|\hat x-\hat y|}{\ep^\beta}\leq C\left(1+ \frac{\ep^{1-\beta}}{\lambda}\right)
\qquad {\rm and }\qquad \frac{|\hat t -\hat s|}{\ep}\leq C\;.
\ee

We claim that
\be\label{barHSub}
\frac{\hat t-\hat s}{\ep} \leq \bar H\left(\frac{\hat x-\hat y}{\ep^\beta}\right)+C\left(\frac{\ep^{1-\beta}}{\lambda}+\lambda\right)\;.
\ee

Indeed for $\alpha,\beta>0$ small, let $(x_\alpha, t_\alpha,y_\alpha,r_\alpha,z_\alpha)$ be a maximum point of $\Psi_1$ given by
$$
\begin{array}{rl}
\Psi_1(x,t,y,r,z)\; = & \displaystyle{  u^\ep(x,t)-\ep v_\lambda\left(y,r,\frac{z-\hat y}{\ep^\beta}\right)
-\frac{|x-\hat y|^2}{2\ep^\beta} -\frac{(t-\hat s)^2}{2\ep}  }\\
&  \displaystyle{ -\frac{1}{2\alpha}\left(|\ep y-x|^2+|z-x|^2+|\ep r-t|^2\right) -\frac{\beta}{2}(|x-\hat x|^2+(t-\hat t)^2) .}
\end{array}
$$

Since $(\hat x,\hat t)$ is the unique maximum point of $\Psi_1(x,t,x/\ep,t/\ep,x)$, we
have that $(x_\alpha, t_\alpha,y_\alpha,r_\alpha,z_\alpha)$ converges to $(\hat x,\hat t, \hat x/\ep, \hat t/\ep, \hat x)$ as $\alpha\to 0$, with
\be\label{Sti2}
\lim_{\alpha\to 0^+} \ \frac{1}{2\alpha}\left(|\ep y_\alpha-x_\alpha|^2+|z_\alpha-x_\alpha|^2 + |\ep r_\alpha - t_\alpha |^2 \right)=0\;,
\ee
while (\ref{LipVlambda}) implies that
\be\label{Sti1}
\frac{|z_\alpha-x_\alpha|}{\alpha}\leq C \frac{\ep^{1-\beta}}{\lambda}\;.
\ee

From the equation satisfied by $u^\ep$ we have
\be\label{alpha1}
\begin{array}{l}
\displaystyle{ \frac{t_\alpha-\hat s}{\ep} + \frac{t_\alpha - \ep r_\alpha }{\alpha}+\beta(t_\alpha-\hat t)}\\
\quad \displaystyle {\leq H\left( \frac{x_\alpha}{\ep}, \frac{t_\alpha}{\ep},
  \frac{x_\alpha-\hat y}{\ep^\beta}+\frac{x_\alpha-\ep y_\alpha}{\alpha}
+\frac{x_\alpha-z_\alpha}{\alpha}+\beta(x_\alpha-\hat x)\right)\;,}
\end{array}
\ee
while from the equation satisfied by $v_\lambda$  we also have
$$
\frac{t_\alpha - \ep r_\alpha}{\alpha} + \lambda v_\lambda\left( y_\alpha, r_\alpha, \frac{z_\alpha-\hat y}{\ep^\beta}\right) \geq
H\left( y_\alpha, r_\alpha, \frac{x_\alpha-\ep y_\alpha}{\alpha}+\frac{z_\alpha-\hat y}{\ep^\beta}\right)\;.
$$
Using the bound on the oscillation of the $\lambda v_\lambda$'s in Lemma \ref{osc}, we get
$$
\bar H\left( \frac{z_\alpha-\hat y}{\ep^\beta}\right) + C\lambda \frac{|z_\alpha-\hat y|}{\ep^\beta}
\geq H\left( y_\alpha, r_\alpha, \frac{x_\alpha-\ep y_\alpha}{\alpha}+\frac{z_\alpha-\hat y}{\ep^\beta}\right) - \frac{t_\alpha - \ep r_\alpha}{\alpha}\;.
$$
Combining the above inequality with (\ref{alpha1}) and using the regularity of $H$ gives
$$
\begin{array}{l}
\displaystyle{  \bar H( \frac{z_\alpha-\hat y}{\ep^\beta}) + C\lambda \frac{|z_\alpha-\hat y|}{\ep^\beta}
\geq \frac{t_\alpha-\hat s}{\ep}+\beta(t_\alpha-\hat t)  }\\
\qquad \qquad \displaystyle{
- C( \frac{|z_\alpha-x_\alpha|}{\ep^\beta}+\frac{|z_\alpha-x_\alpha|}{\alpha}
+ \beta|x_\alpha-\hat x| )} \\
\qquad \qquad \displaystyle{
- C(\left|\frac{x_\alpha}{\ep}-y_\alpha\right| + \left|\frac{t_\alpha}{\ep}-r_\alpha\right| ) ( \frac{|x_\alpha-\ep y_\alpha|}{\alpha}
+\frac{|z_\alpha- \hat y|}{\ep^\beta}).   }
\end{array}
$$
Using (\ref{Sti2}) and (\ref{Sti1}) we now let $\alpha\to 0$  to get
$$
\bar H\left( \frac{\hat x-\hat y}{\ep^\beta}\right) + C\lambda \frac{|\hat x-\hat y|}{\ep^\beta}
\geq \frac{\hat t-\hat s}{\ep} - C\frac{\ep^{1-\beta}}{\lambda}\;.
$$

Recalling (\ref{Sti0}) we obtain
$$
\bar H\left( \frac{\hat x-\hat y}{\ep^\beta}\right) + C \left(\lambda+ \ep^{1-\beta}+\frac{\ep^{1-\beta}}{\lambda}\right)
\geq \frac{\hat t-\hat s}{\ep},
$$

and, since $\lambda\in (0,1)$,   we finally get (\ref{barHSub}). \\

We now show that
\be\label{barHSup}
\frac{\hat t-\hat s}{\ep} -\delta \geq \bar H\left(\frac{\hat x-\hat y}{\ep^\beta}\right) -C\frac{\ep^{1-\beta}}{\lambda}\;.
\ee

To this end, for $\alpha,\beta>0$, we consider
$$
\begin{array}{rl}
\Psi_2(y,s,z)\; = & \displaystyle{ \bar u(y,s)+\ep v_\lambda\left(\frac{\hat x}{\ep},\frac{\hat t}{\ep},\frac{\hat x-z}{\ep^\beta}\right)
+\frac{|\hat x-y|^2}{2\ep^\beta} +\frac{(\hat t-s)^2}{2\ep}+\delta s }\\
&\qquad\qquad \displaystyle{  +\frac{ |z-y|^2}{2\alpha}
+\frac{\beta}{2}(|y-\hat y|^2+|s-\hat s|^2),  }
\end{array}
$$
which has a minimum at some point $(y_\alpha,s_\alpha, z_\alpha)$. Using (\ref{LipVlambda}) and the fact that
 $(\hat y, \hat s)$ is the unique minimum of the map $(y,s)\to \Psi_2(y,s,y)$, we find that 
 $(y_\alpha,s_\alpha, z_\alpha)$ converges to $(\hat y, \hat s, \hat y)$ as $\alpha\to 0$, with
$$
\frac{ |z_\alpha-y_\alpha|}{\alpha} \leq C\frac{\ep^{1-\beta}}{\lambda}\qquad {\rm and }\qquad
\lim_{\alpha\to0} \frac{ |z_\alpha-y_\alpha|^2}{\alpha}=0\;.
$$

Since $\bar u$ solves \eqref{eq:hom} we have
$$
\frac{\hat t-s_\alpha}{\ep}-\delta \geq \bar H (\frac{\hat x-y_\alpha}{\ep^\beta}-\frac{y_\alpha-z_\alpha}{\alpha}
-\beta(y_\alpha-\hat y))\;,
$$
and, in view of the Lipschitz continuity of $\bar H$, letting $\alpha\to0^+$ gives
(\ref{barHSup}). \\

Combining (\ref{barHSub}) and (\ref{barHSup}) we obtain
$$
\delta-C\left( \frac{\ep^{1-\beta}}{\lambda}+\lambda\right)\leq 0\;,
$$
which for  $\lambda= \ep^{(1-\beta)/2}$ and $\delta=3C \ep^{(1-\beta)/2}$  gives a contradiction. So either $\hat t=0$
or $\hat s=0$.\\

Next we estimate $\max_{x,y,s} \Phi(x,0,y,s)$ and $\max_{x,t,y} \Phi(x,t,y,0)$.
We have:
$$
\begin{array}{l}
\displaystyle{ \max_{x,y,s} \Phi(x,0,y,s)}\\
\qquad \qquad \leq \; \displaystyle{ \max_{x,y,s} \left\{u_0(x)-u_0(y)+Cs
+C\frac{\ep^{1-\beta}}{\lambda}|x-y|
-\frac{|x-y|^2}{2\ep^\beta} -\frac{s^2}{2\ep}-\delta s \right\}  }\\
\qquad \qquad \leq \; \displaystyle{ \max_{x,y} \left\{C|y-x|+C \ep
-\frac{|x-y|^2}{2\ep^\beta}  \right\}  }\\
\qquad \qquad \leq \; \displaystyle{ C\left(  \ep^\beta   +\ep \right) \; \leq C\epsilon^\beta , }\\
\end{array}
$$
and, similarly,
$$\displaystyle{
\max_{x,y,s} \Phi(x,t,y,0) \leq C\ \ep^\beta }.$$

Therefore, for any $(x,t)\in \R^N\times [0,T]$, we have
$$
u^\ep(x,t)-\bar u(x,t)\;  \leq \; \delta t +C\epsilon^\beta \; \leq\; C T\ep^{(1-\beta)/2}+C\ep^\beta\;.
$$

Choosing $\beta=1/3$ finally gives
$$
u^\ep(x,t)-u(x,t) \leq  C\ep^{1/3}\;.
$$

This reverse inequality can be obtained in the same way.
\QED

%%%%%%%%%%%%%%%%%%%%%%%%%%%%%%%%%%%%%%%%%%%%%%%%%%
%%%%%%%%%%%%%%%%%%%%%%%%%%%%%%%%%%%%%%%%%%%%%%%%%%
\section{Enhancement of speed}\label{sec:enh}

To prove Theorem \ref{EnhenVx} and Theorem \ref{EnhenVxt} we need three  results
which we formulate as Lemmas next but present their proof at the end of the section.

We begin with

\begin{Lemma}\label{lossgen} Let $V$ and $\bar H$ be as in Theorem \ref{theo:mainperiodic}.
Then, for any  $\Z-$periodic $c\in {\mathcal C}^1(\R;\R^N)$,
the averaged Hamiltonian associated to $V-c$
is  $\bar H-\lg \int_0^1 c(s)ds,\cdot\rg $.
\end{Lemma}

The second is

\begin{Lemma}\label{constant} The divergence zero condition (\ref{divzPV0}) is equivalent to the fact that, for any fixed time $t$,
the map $x\to z(x,t)+\lg P,x\rg$ is constant along the flow of the ode
$X'(s)=V(X(s),t)$.
\end{Lemma}

To state the final result we recall the notion of an $\ep$-mollifier. To this end, let $\phi \in {\mathcal C}_c^{\infty}(\R^{N+1};[0,1])$
be such that $\phi (0)=1$ and $\int_{\R^{N+1}}\phi=1$ and define the $\ep$-mollifier $\phi_\ep$ by $\phi_\ep = \ep^{-(N+1)} \phi (x/\ep, t/\ep)$. Then
$\int_{\R^{N+1}}\phi_\ep=1$ and, for any $f \in L_{loc}^1(\R^{N+1})$, $f_\ep = f * \phi_\ep$ is a smooth approximation of $f$.

For the rest of the section we assume that, for each $t\in\R$,
\be\label{eq:extra}
{\rm div}_xV(\cdot,t)=0 \quad \text{ and } \quad \int_{Q_1}V(x,t)dx=0;
\ee
recall that the average zero condition is actually not a restriction, since we can always
replace $V$ by $V-\int_{Q_1}V(x,t)dx$ a fact, which, in view of Lemma $3.1$, simply adds a translation
to the effective Hamiltonian.

The final preliminary result is

\begin{Lemma}\label{lem:SuperCorrector} Assume  \eqref{eq:extra} and
%that, for all $t\in\R$, ${\rm div}_xV(\cdot,t)=0$ and $\int_{Q_1}V(x,t)dx=0$ and
fix $P\in\R^N$. There exits a bounded, $\Z^{N+1}$-periodic
$z\in BV_{loc}(\R^N\times \R)$ such that, for all $\ep>0$
%and
%all compactly supported in $B_\ep(0)$ nonnegative kernels $\phi_\ep\in {\mathcal C}^\infty(\R^{N+1})$ such that $\int_{\R^{N+1}}\phi_\ep=1$,
the smooth functions
$z_\ep=z*\phi_\ep$'s satisfy in $\R^{N+1}$
\be\label{IneqSuperCorrector}
\begin{array}{l}
\ds{ z_{\ep,t}+ \bar H(P) \geq  |Dz_{\ep}+P| }\\
\quad \ds{ +
\int_{\R^{N+1}}(z(y,s)+\lg P,y\rg)\lg D\phi_\ep(\cdot-y,\cdot-s), V(y,s)\rg dyds \;.}
\end{array}
\ee
\end{Lemma}

%For the rest of the section we assume that, for each $t\in\R$,
%\be\label{eq:extra}
%{\rm div}_xV(\cdot,t)=0 \quad \text{ and } \quad \int_{Q_1}V(x,t)dx=0;
%\ee
%recall that the average zero condition is actually not a restriction, since we can always
%replace $V$ by $V-\int_{Q_1}V(x,t)dx$ a fact, which, in view of Lemma $3.1$, simply induces a translation
%to the effective Hamiltonian.

We proceed with Theorem \ref {EnhenVx}, which is a straightforward consequence of  Theorem \ref{EnhenVxt}. We have:

\noindent {\bf Proof of Theorem \ref{EnhenVx}:} If $\bar H(P)=|P|$, then
let 
$$z(x,t)=\hat z(\frac{\lg P,x\rg}{|P|}+t)\quad \text{ and } \quad z_1(x)=\int_0^1 z(x,s)ds,$$
where $\hat z$ is given by Theorem \ref{EnhenVxt}.

Since $V$ is independent of $x$, $z_1$ also
satisfies (\ref{divzPV0}) in the sense of distributions, and, since $z$ is periodic $z_1$ is actually a constant.
In this case (\ref{divzPV0}) reduced to $\lg V(x), P\rg =0$ for all $x\in \R^N$.

When $N=2$ and $\lg V \rg=0$ there exists a $\Z^2$-periodic stream function $E:\R^2\to \R$ such that
$V=(-\frac{\partial E}{\partial x_2},\frac{\partial E}{\partial x_1})$. 
In this case, if $P=(P_1,P_2)$ and $q=(-P_1, P_2)$, $\lg V, P\rg=0$ in $\R^2$ becomes
$$
0= \lg DE(x),q \rg \;.
$$

If $q$ is an irrational direction, then the map $t\to x+tq$ is dense in $\R^2/Z^2$. Since $E$ is constant along this
trajectory, $E$ is constant
and therefore $V$ is identically equal to $0$. Otherwise, $t\to x+tq$ is periodic and $E$ has to be constant along this trajectory.
This means that
$E= \tilde E(\lg \cdot , P\rg)$ for some smooth periodic map $\tilde E:\R\to \R$, and, hence, $V$ is a shear advection.

Conversely, if $\lg V, P\rg=0$ in $\R^N$, then $v_\lambda=|P|/\lambda$ is the unique solution to (\ref{eq:vlambda}),
and $\lambda v_\lambda = |P|$ clearly converges uniformly to $\bar H(P)=|P|$.
\QED

We turn now to the

\noindent {\bf Proof of Theorem \ref{EnhenVxt}:} Assume that, for some $P\in \R^N\backslash \{0\}$, $\bar H(P)=P$  and let
$z$ be given by Lemma \ref{lem:SuperCorrector}. We first  prove that $z$ is a function of only $\lg P, x\rg$
and $t$. More precisely, we claim that there exits $\tilde z\in BV_{loc}(\R^2)$ such that
\be\label{z=zpxt}
\mbox{\rm $z=\tilde z(\lg P, x\rg,t)$ and $ \tilde z_s(s,t)\geq -1$ in the sense of distributions. }
\ee

To this end, let $\phi_\ep$ and $z_\ep=z*\phi_\ep$
be as in Lemma \ref{lem:SuperCorrector}; note that $z_\ep$ is $\Z^{N+1}$-periodic. Then, for all $(x,t)\in \R^N\times\R$,
\be \label{eq:triou}
\int_{Q_1} \int_{\R^{N+1}}(z(y,s)+\lg P,y\rg)\lg D\phi_\ep(x-y,t-s), V(y,s)\rg dyds\ dx =0\;.
\ee

Indeed the periodicity of $z$ and $V$ gives,
$$
\begin{array}{l}
\displaystyle{  \int_{Q_1} \int_{\R^{N+1}}(z(y,s)+\lg P,y\rg)\lg D\phi_\ep(x-y,t-s), V(y,s)\rg dyds \ dx }\\
\quad \displaystyle{  = \int_{\partial Q_1} \int_{\R^{N+1}}(z(y,s)+\lg P,y\rg) \phi_\ep(x-y,t-s) \lg V(y,s),\nu_x \rg dyds \ d{\mathcal H}^{N-1}(x)  }\\
\quad \displaystyle{  =  \int_{\R^{N+1}}\int_{\partial Q_1} (z(x-y,t-s)+\lg P,x-y\rg) \phi_\ep(y,s) \lg V(x-y,t-s),\nu_x \rg d{\mathcal H}^{N-1}(x)\  dyds  }\\
\quad \displaystyle{  =  \int_{\R^{N+1}}\int_{\partial Q_1} \lg P,x-y\rg \phi_\ep(y,s) \lg V(x-y,t-s),\nu_x \rg d{\mathcal H}^{N-1}(x)\ dy ds. }\\
\end{array}
$$
while \eqref{eq:extra} yields

$$
\begin{array}{l}
\displaystyle{  \int_{\R^{N+1}}\int_{\partial Q_1} \lg P,x-y\rg \phi_\ep(y,s) \lg V(x-y,t-s),\nu_x \rg d{\mathcal H}^{N-1}(x)\ dy ds }\\
\qquad \displaystyle{  =   \int_{\R^{N+1}}\int_{ Q_1} {\rm div}_x \left(\lg P, x-y\rg \phi_\ep(y,s) V(x-y,t-s)\right)  dx\ dy ds }\\
\qquad \displaystyle{  =   \int_{\R^{N+1}}\int_{ Q_1} \lg P, V(x-y,t-s)\rg \phi_\ep(y,s)  dx \ dy ds\;=\; 0 \;,}\\
\end{array}
$$
and, hence, \eqref{eq:triou} holds.

Next we
integrate (\ref{IneqSuperCorrector}) over $Q_1\times (0,1)$. Using \eqref{eq:triou} and the periodicity of $z_\ep(\cdot,t)$ we get
$$
|P|= \bar H(P) \geq \int_0^1\int_{Q_1} |D z_\ep(y,t)+P|dydt  \geq \int_0^1\left| \int_{Q_1}(D z_\ep(y,t)+P)dy\right|dt= |P|\;.
$$

It follows that, for all $(x,t)\in \R^{N+1}$, there exists $\theta(x,t)\geq -1$ such that $Dz_\ep(x)=\theta(x,t)P$.
Thus $z_\ep$ is of the form $z_\ep(x,t)=\tilde z_\ep(\lg x,P \rg,t)$, with $\tilde z_\ep:\R^2\to \R$ satisfying
$\tilde z_{\ep,s}(s,t)\geq -1$.

Passing to the limit $\ep\to0$, we also find that
$z= \tilde z(\lg x,P\rg,t)$ for some map $\tilde z\in BV_{loc}(\R^2,\R)$ satisfying
$\tilde z_s(s,t)\geq -1$ in the sense of distributions. Whence (\ref{z=zpxt}) holds. \\

Next we claim that $\tilde z$ satisfies, in the sense of distributions,
\be\label{transpz}
\tilde z_t(s,t)=\tilde z_s(s,t)|P|\qquad  {\rm in }\; \R^2 ,
\ee
and that (\ref{divzPV0}) holds. Note that this proves the ``if'' part, since (\ref{transpz}) implies the existence of a
 map $\hat z\in BV_{loc}(\R,\R)$ such that
$$
\tilde z(s,t)= \hat z (|P|^{-1}s + t)\;.
$$

Moreover we have $\hat z'(s) \geq -|P|$ in the sense of distributions because $ \tilde z_s(s,t)\geq -1$ in the same sense.
Finally, $z(x,t)=\hat z\left(\frac{\lg P, z\rg}{|P|}+t\right)$ is periodic in space and time.\\

We continue with the proofs of (\ref{transpz}) and (\ref{divzPV0}). If $z$ is constant, then (\ref{transpz}) is obvious and
(\ref{divzPV0}) just follows from (\ref{IneqSuperCorrector}) when $\ep\to0$.

Next we assume that $z$ is not constant. In this case $z= \tilde z(\lg \cdot,P\rg,t)$ is $\Z^{N+1}$-periodic and not constant.
Therefore $P$ has to be a rational direction. So, up to a rational change of coordinates, we may assume without loss of generality
that $P=\theta e_1$ for some $\theta>0$, while $V$ is still $\Z^{N+1}-$periodic. 

Using the notation $x=(x_1,x')$
for each vector of $\R^N$ with $x_1\in \R$ and $x'\in \R^{N-1}$, for a fixed $(x_1,t)\in \R^2$
we integrate (\ref{IneqSuperCorrector}) over  the cube $\{x_1\}\times Q_1'\times\{t\}$, where
$Q_1'=(-1/2,1/2)^{N-1}$, and obtain
\be\label{TATAT}
\begin{array}{l}
\ds{\tilde z_{\ep,t}(\lg P,x\rg,t)+ |P| \; \geq \; (\tilde z_{\ep,s}(\lg P,x\rg,t) +1)|P| + }\\
 \ds{ \int_{Q_1'}\int_{\R^{N+1}}(\tilde z(\theta(x_1-y_1),t-s)+\theta(x_1-y_1))\lg D\phi_\ep(y,s), V(x-y,t-s)\rg dyds\ dx'\;.}
\end{array}
\ee

It turns out that the last integral in the right-hand side of the above inequality vanishes. Indeed
$$
\begin{array}{l}
\ds{  \int_{Q_1'}
\int_{\R^{N-1}}\lg D\phi_\ep(y,s), V(x-y,t-s)\rg dy' dx' \; = \;} \\
\ds{
\int_{\R^{N-1}} \frac{\partial\phi_\ep}{\partial x_1}(y,s)\int_{Q_1'} V_1(x-y,t-s) dx' dy'}\\
\ds{ +
\sum_{j=2}^N \int_{Q_1'}
\int_{\R^{N-1}}\phi_\ep(y,s) \frac{\partial  V_j}{\partial x_j}(x-y,t-s)  dy' dx'  \;. }
\end{array}
$$

The periodicity of $V$ yields, for any $j = 2, \dots, N$,
$$
\int_{\R^{N-1}}\phi_\ep(y,s) \int_{Q_1'} \frac{\partial  V_j}{\partial x_j}(x-y,t-s)   dx'dy'
=
0\;,
$$
while the divergence free condition and, again, the periodicity give
$$
\frac{\partial}{\partial x_1} \int_{Q_1'}
V_1(x-y,t-s)  dx'
=
\int_{Q_1'} \frac{\partial V_1}{\partial x_1}(x-y,t-s)  dx'
=-\sum_{j=2}^N \int_{Q_1'}\frac{\partial V_1}{\partial x_j}(x-y,t-s)  dx'=0\ .
$$

On the other hand,
$$
\int_{-1/2}^{1/2} \int_{Q_1'}
V_1(x-y,t-s)  dx' \,dx_1=\int_{Q_1} V_1(x,t-s)dx =0,
$$
and, hence, for all $x_1\in \R$,
$$
\int_{Q_1'}
V_1(x-y,t-s)  dx'
=0 \;.
$$

Therefore
$$
\int_{Q_1'}
\int_{\R^{N-1}}\lg D\phi_\ep(y,s), V(x-y,t-s)\rg dy' dx'=0\;,
$$
which, going back to (\ref{TATAT}), proves that
$$
\tilde z_{\ep,t}(\lg P,x\rg,t) \geq  \tilde z_{\ep,s}(\lg P,x\rg,t)|P| \quad \text { in } \quad \R^{N+1}\;.
$$

Since $\tilde z_\ep$ is periodic, integrating the above inequality over $(-1/2,1/2)\times (0,1)$ shows that in fact it must be an equality.
Letting $\ep\to0$ then gives (\ref{transpz}).

To prove (\ref{divzPV0}), we first combine  (\ref{IneqSuperCorrector}), (\ref{z=zpxt}) and (\ref{transpz}) to get, 
for all $(x,t)\in \R^{N+1}$,
$$
0\geq \int_{\R^{N+1}}(z(y,s)+\lg P,y\rg)\lg D\phi_\ep(x-y,t-s), V(y,s)\rg dyds \;.
$$

Averaging over the cube $Q_1$ we see that, as a matter of fact, equality must hold for all $(x,t)\in \R^{N+1}$.

Integrating the resulting equality against any compactly supported smooth function $\psi:\R^N\to \R$ we get
$$
\int_{\R^{N}}\int_{\R^{N+1}}\psi(x)(z(y,s)+\lg P,y\rg)\lg D\phi_\ep(x-y,t-s), V(y,s)\rg dyds \ dx\;=\; 0 \;,
$$
and after integrating again by parts and letting $\ep\to0$, we obtain, for all $t\in\R$,
$$
\int_{\R^{N}}\lg D\psi(x), (z(x,t)+\lg P,x\rg)V(x,t)\rg dx\;=\; 0 \;,
$$
which is exactly (\ref{divzPV0}).

To prove the ``only if" part let $\hat z$ be as claimed.
Since $\hat z\in BV_{loc}(\R,\R)$, we may assume without loss of generality
that $\hat z$, and, hence, $z(x,t)=\hat z(\frac{\lg P,x\rg}{|P|}+t)$ are lower semi-continuous.
% so that $z(x,t)=\hat z(\frac{\lg p,x\rg}{|P|}+t)$ is also
%lower  semi continuous.

We show next that $z$ satisfies, in the viscosity sense, 
\be\label{CorrEq2}
\partial_t z+|P|\geq |Dz+P|+ \lg V, Dz+P\rg\qquad {\rm in}\; \R^{N+1}\;.
\ee

To this end, let $\phi$ be a smooth test function such that $z\geq \phi$ with equality at $(\bar x,\bar t)$.
It follows from equality $z(x,t)=\hat z(\frac{\lg P,x\rg}{|P|}+t)$
that $D\phi(\bar x,\bar t)=\theta P/|P|$ where $\theta=\phi_t(\bar x,\bar t)$.

Since
$\hat z'\geq -|P|$ in the sense of distribution, it follows that $\theta \geq -|P|$. Moreover recalling (\ref{divzPV0}) and Lemma 
\ref{constant}, we have, for a fixed $t$, that the function $x\to V(x,t)+\lg P,x\rg$ is constant under the flow of the ode
$X'(s)=V(X(s),t)$.

Let now $X$ be a solution with $X(0)=\bar x$
and $t=\bar t$. Then, for any $s\in \R$,
$$
\phi(\bar x,\bar t)+\lg P, \bar x\rg = z(\bar x,\bar t)+\lg P,\bar x\rg
= z(X(s),\bar t)+\lg P,X(s)\rg
\geq  \phi(X(s),\bar t)+\lg P,X(s)\rg \;,
$$
and, therefore,
$$
0=\frac{d}{ds}_{|_{s=0}}\left[ \phi(X(s),\bar t)+\lg X(s),P\rg\right] =\lg D\phi (\bar x,\bar t)+P, V(\bar x,\bar t)\rg\;.
$$

Combining the above relations gives
$$
\phi_t(\bar x,\bar t)+|P|= |D\phi(\bar x,\bar t)+P|+ \lg V(\bar x,\bar t), D\phi(\bar x,\bar t)+P\rg\;,
$$
i.e., $z$ is a super-solution of (\ref{CorrEq2}).

%For any $\lambda >0$, let $v_\lambda$ be the solution of \eqref{eq:vlambda}.
%\be\label{eq:cell}
%\partial_t v_\lambda + \lambda v_\lambda = |Dv_\lambda+P|+\lg V(x,t), Dv_\lambda+P \rg \qquad {\rm in }\; \R^{N+1}\;.
%\ee

Then it is easy to check that $\tilde z (x,t)= z(x,t)+\frac{|P|}{\lambda}+\|z\|_\infty$ is a super-solution of \eqref{eq:vlambda}. It then
follows from the comparison principle that
$$v_\lambda \leq z+ {\lambda}^{-1} |P|+\|z\|_\infty\;.
$$

Recalling that $\lambda v_\lambda$ converges uniformly, as $\lambda\to0$, to $\bar H(P)$, we obtain
$\bar H(P)\leq |P|$. Since the reverse inequality always holds, the proof of the implication is complete.
\QED

We continue with the proofs of the lemmas. We begin with the

\noindent{\bf Proof of Lemma \ref{lossgen}:}
Fix $\lambda >0$ and $P\in \R^N$, let $v_\lambda$ be the solution of \eqref{eq:vlambda}
%$$
%v_{\lambda,t} + \lambda v_\lambda = |Dv_\lambda+P|+\lg V, Dv_\lambda+P \rg \qquad {\rm in }\; \R^{N+1},
%$$
and us recall that $\lambda v_\lambda$ converges uniformly, as $\lambda \to 0$, to $\bar H(P)$.

Set $\bar c=\int_0^1 c(s)ds$ and consider $w_\lambda \in BUC(\R^{N+1})$ given by
$$
w_\lambda(x,t)=v_\lambda(x-\int_0^t c(s)ds, t)-\lg \int_0^t(c(s)-\bar c)ds+ \frac{\bar c}{\lambda},P\rg +2\|c\|_\infty|P|\;.
$$

It follows that $w_\lambda$ is a super-solution of %\eqref{eq:lossgen}
\be\label{eq:lossgen}
z_{\lambda,t} + \lambda z_\lambda =  |Dz_\lambda+P|+\lg V-c, Dz_\lambda+P \rg \qquad {\rm in }\; \R^{N+1}.
\ee

We only present a formal proof, which can be easily justified using viscosity solution arguments.
To this end observe that it is immediate from the definition of $w_\lambda$ that
$$
w_{\lambda,t} +\lambda w_\lambda\geq  -\lg Dv_\lambda , c(t)\rg + v_{\lambda,t}- \lg c-\bar c,P\rg
+ \lambda v_\lambda- \lg \bar c,P\rg ,
$$
while
$$
|Dw_\lambda+P|+\lg V -c, Dw_\lambda+P \rg
=
|Dv_\lambda+P|+\lg V - c, Dv_\lambda+P \rg\; ;
$$
the comparison principle now gives  $w_\lambda\geq z_\lambda$ where
$z_\lambda$ is the solution of \eqref{eq:lossgen}.

Since the $\lambda z_\lambda$'s and the  $\lambda w_\lambda$'s converge uniformly, as $\lambda\to 0$,
to the averaged Hamiltonian $\bar H_c(P)$ associated to $V-c$ and
to $\bar H(P)-\lg \bar c, P\rg$ respectively, we get $\bar H_c(P)\leq \bar H(P)-\lg \bar c, P\rg$.
The opposite inequality is proved similarly by considering $-c$ instead of $c$.
\QED

We continue with the

\indent {\bf Proof of Lemma \ref{constant} : } Let $X_x(\cdot)$ be the solution of $X'(s)=V(X(s),t)$ with initial condition
$x$ at time $s=0$. Then, for any $h\in \R$ and $\psi:C_c^{\infty}(\R^N)$ the divergence zero
property of $V$ yields
$$
\int_{\R^N}\psi(x) (z(X_x(h),t)+\lg P, X_x(h)\rg) dx=\int_{\R^N}\psi(X_x(-h)) (z(x,t)+\lg P, x\rg) dx\;.
$$

Therefore, in view of (\ref{divzPV0}),
$$
\begin{array}{l}
\ds{ \frac{d}{dh}_{|_{h=0}}\int_{\R^N}\psi(x) (z(X_x(h),t)+\lg P, X_x(h)\rg) dx}\\
\qquad\qquad \ds{ =\int_{\R^N}\lg D\psi(x),V(x,t)\rg (z(x,t)+\lg P, x\rg) dx=0\;,}
\end{array}
$$

Applying this last equality to the test function $\psi\circ X_x(-s)$ we get
$$
\begin{array}{l}
\ds{ \frac{d}{dh}_{|_{h=0}}\int_{\R^N}\psi(x) (z(X_x(s+h),t)+\lg P, X_x(s+h)\rg) dx}\\
\qquad\qquad \ds{ =\frac{d}{dh}_{|_{h=0}}\int_{\R^N}\psi(X_x(-s)) (z(X_x(h),t)+\lg P, X_x(h)\rg) dx= 0\;. }
\end{array}
$$

Hence $\ds{ \int_{\R^N}\psi(x) (z(X^x(s),t)+\lg P, X^x(s)\rg) dx  }$ is constant in time, which means that
$z(\cdot,t)$ is constant along the flow.
\QED

We conclude with the

\noindent {\bf Proof of Lemma \ref{lem:SuperCorrector} : } For $\lambda >0$, let $v_\lambda$ be the solution of
of \eqref{eq:vlambda} 
%$$
%\partial_t v_\lambda + \lambda v_\lambda = |Dv_\lambda+P|+\lg V(x,t), Dv_\lambda+P \rg \qquad {\rm in }\; \R^{N+1}
%$$
and set
$$z_\lambda(x,t)=v_\lambda(x,t)-v_\lambda(0,0).$$

It follows from Lemma \ref{osc} and Corollary \ref{cor:barH} that, for some $C>0$ independent of $\lambda$,
$$\|z_\lambda\|_\infty \leq C|P|.$$

Next we show that the
 $z_{\lambda}$'s  are also bounded in $BV_{loc}(\R^{N+1})$. Indeed, for $\alpha>0$, consider
the $\Z^{N+1}$-periodic solution $v_{\lambda,\alpha}$ to
\be\label{vlambdaalpha}
v_{\lambda,\alpha,t}+
\lambda v_{\lambda ,\alpha}=\alpha \Delta v_{\lambda, \alpha}+|Dv_{\lambda ,\alpha}+P|+\lg V, Dv_{\lambda ,\alpha}+P\rg \qquad {\rm in}\; \R^{N+1}\;,
\ee
which is at least in $C^{1,1}$ and, moreover, converges uniformly, as $\alpha\to0$, to
$v_\lambda$.

Integrating (\ref{vlambdaalpha}) over a cylinder of the form $Q_R\times (-R,R)$ for some positive integer $R$, we obtain,
using the periodicity, that
$$
\begin{array}{l}
\ds{ \int_{Q_R\times (-R,R)} \lambda v_{\lambda ,\alpha} \geq \int_{Q_R\times (-R,R)} |Dv_{\lambda ,\alpha}+P| }\\
\qquad\qquad \ds{  +
\int_{Q_R\times (-R,R)} \lg V(x,t), Dv_{\lambda ,\alpha}(x,t)+P\rg dxdt \;. }
\end{array}
$$

Since, in view of \eqref{eq:extra},
$$
\int_{Q_R} \lg V(x,t), Dv_{\lambda ,\alpha}(x,t)+P\rg  dx =0,
$$
it follows that
\be \label{BVboundx}
\int_{Q_R\times (-R,R)} |Dv_{\lambda ,\alpha}+P| \leq 2R^{N+1}\|\lambda v_{\lambda ,\alpha}\|_\infty\;.
\ee

Let $\Phi\in C_c^{\infty}(Q_R\times (-R,R);\R^{N+1})$ be such that
$|\Phi(x,t)|\leq 1$ for all $(x,t)\in \R^{N+1}$. From (\ref{vlambdaalpha}) and (\ref{BVboundx}) we get
$$
\begin{array}{l}
\ds{ \int_{Q_R\times (-R,R)} v_{\lambda,\alpha} {\rm div}_{x,t}\Phi dxdt
=-\int_{Q_R\times (-R,R)} [v_{\lambda,\alpha,t} \Phi + \lg D v_{\lambda,\alpha}, \Phi\rg ]dxdt =} \\
\ds{-\int_{Q_R\times (-R,R)}
[(\lambda v_{\lambda ,\alpha}-\alpha \Delta v_{\lambda, \alpha}-|Dv_{\lambda ,\alpha}+P|-\lg V , Dv_{\lambda ,\alpha}+P\rg ) \Phi
+ \lg D v_{\lambda,\alpha}, \Phi \rg]dxdt } \\
\ds{ \leq C_R(P) \| \lambda v_{\lambda ,\alpha}\|_\infty +\alpha \int_{Q_R\times (-R,R)} v_{\lambda, \alpha}\Delta \Phi dxdt },
\end{array}
$$
where $C_R(P)$ depends only on $N$, $R$, $\|V\|_\infty$ and $P$.

Letting $\alpha\to 0$ yields
\be\label{eq:extra1}
\int_{Q_R\times (-R,R)} v_{\lambda} {\rm div}_{x,t}\Phi dxdt \leq C_R(P) \| \lambda v_{\lambda }\|_\infty
\leq C_R(P)(\bar H(P) +C\lambda) \;,
\ee
which in turn implies, 
in view of the assumptions on $\Phi$, that
the $z_\lambda$'s are bounded in $L^\infty$ and in $BV_{loc}$ uniformly with respect to
$\lambda$. Hence the $z_\lambda$'s converge, up to a subsequence and in $L^1_{loc}$, to some $\Z^{N+1}$- periodic
$z\in BV_{loc}(\R^{N+1})$.

Let $\ep>0$, $\phi_\ep$ and $z_{\ep}$ as in statement of the lemma,
set $z_{\lambda, \ep}= \phi_\ep*z_{\lambda }$, and
$z_{\lambda, \alpha,\ep}= \phi_\ep*z_{\lambda ,\alpha}$ and fix  $(x,t)\in \R^{N+1}$.

It follows from (\ref{vlambdaalpha}) that
$$
\begin{array}{l}
\ds{ z_{\lambda ,\alpha,\ep,t}+ \lambda z_{\lambda ,\alpha,\ep}+ \lambda v_{\lambda ,\alpha}(0,0)
\geq \alpha \Delta z_{\lambda ,\alpha,\ep}+|Dz_{\lambda ,\alpha,\ep}+P|}\\
\qquad \ds{ +
\int_{\R^{N+1}}\phi_\ep(x-y,t-s)\lg V(y,s), Dz_{\lambda ,\alpha}(y,s)+P\rg dyds\;.}
\end{array}
$$

Integrating by parts and using \eqref{eq:extra} we find
$$
\begin{array}{l}
\ds{  z_{\lambda ,\alpha,\ep, t}+ \lambda z_{\lambda ,\alpha,\ep}+ \lambda v_{\lambda ,\alpha}(0,0)
\geq \alpha \Delta z_{\lambda ,\alpha,\ep}+|Dz_{\lambda ,\alpha,\ep}+P|}\\
\qquad \ds{  +\int_{\R^{N+1}}(z_{\lambda,\alpha}(y,s)+\lg P,y\rg)\lg D\phi_\ep(x-y,t-s), V(y,s)\rg dyds\;. }
\end{array}
$$

Letting first $\alpha\to 0$ and then $\lambda\to0$ gives (\ref{IneqSuperCorrector}).
\QED

%%%%%%%%%%%%%%%%%%%%%%%%%%%%%%%%%%%%%%%%%%%%%%%%%%
%%%%%%%%%%%%%%%%%%%%%%%%%%%%%%%%%%%%%%%%%%%%%%%%%%
\section{Convergence to the Wulff shape} \label{sec:Wulff}

%Throughout this section we assume that, for all $x\in\R^N$ and all $t\in\R$,
%\be\label{eq:ass}
%{\rm div}_x V(x,t)=0 \quad \text{ and } \quad  \int_{Q_1}V(x,t)dx=0.
%\ee

We begin recalling 
some important facts from the theory of front propagation. The first is (see, for instance, \cite{bss93}), that
the family of sets $(K(t))_{t\geq 0}$ is  independent of the choice of $u_0$ as long as $K_0=\{x\in\R^n : u_0(x) \geq0 \}$.
The second (again see \cite{bss93} and the references therein), which we will use repetitively in the sequel, is the following superposition principle of the geometric
flow. If $(K^\theta_0)_{\theta\in \Theta}$ is an arbitrary family of non-empty closed
subsets of $\R^N$, with corresponding solution $(K^\theta(t))_{t\geq 0}$, then the solution starting from $\overline{\bigcup_{\theta\in \Theta}K^\theta_0}$
is given by $(\overline{\bigcup_{\theta\in \Theta}K^\theta(t)})_{t\geq 0}$.
This can be seen either by using the control representation of the geometric flow or the stability and comparison properties of viscosity solutions. A consequence is the well known inclusion principle.
If $K_0\subset K_0'$ are two non-empty  closed subsets of $\R^N$, then the corresponding solutions
$(K(t))_{t\geq 0}$ and $(K'(t))_{t\geq 0}$ satisfy $K(t)\subset K'(t $) for all $t\geq 0$. \\

We also remark that the convergence, as $t\to\infty$, of the level sets of geometric equations without spatio-temporal inhomogeneities was
considered in \cite{IPS}. The results of \cite{IPS} do not,  however, apply  to the problem at hand.

The proof of Theorem \ref{theo:wulff} is long. We formulate two important steps as separate lemmas which we prove at the end of the section.

The first is about some ``controllability'' estimates. We have:

\begin{Lemma}\label{envahisseur} Assume \eqref{eq:extra}. There exist a positive integer $n_0$ and $T>0$ such that, for all $x\in \R^N$, the solution
$(\hat K(t))_{t\geq0}$ of the front propagation problem starting from the set $x+\{k\in \Z^N\;:\; |k|\leq n_0\}$ contains $Q_1(x)$ at time $T$.
\end{Lemma}

The second is about some growth property for fronts.

\begin{Lemma}\label{lem:keygrowth} There exits a positive integer $\bar R$ and positive constants  $r$ and $T_1$ such that,
 if the initial compact set $K_0\subset \R^N$ contains a set of the form $Q_{\bar R}(k)$
for some $k\in \Z^N$, then the solution $(K(t))_{t\geq 0}$ of the front propagation problem starting from $K_0$ satisfies, for all $t\geq0$,
$$
Q_{rt}(k)\subset K(t+T_1)\qquad \forall t\geq 0\;.
$$
\end{Lemma}

We continue with the

\noindent{\bf Proof of Theorem \ref{theo:wulff} :}
We first show (\ref{InclWulff1}). It is well known (see \cite{bss93})
 that
the characteristic function ${\bf 1}_{K(t)}$ of $K(t)$ is a solution to the geometric equation
\be\label{eq:geoWulff}
\partial_t u=|Du|+\lg V, Du\rg\qquad {\rm in }\; \R^N\times (0,+\infty)\; .
\ee

Fix next a  direction $\nu\in \R^N$ with $|\nu|=1$. For $\lambda>0$, let $v_\lambda$ be the solution of
$$
v_{\lambda,t} + \lambda v_\lambda = |Dv_\lambda+\nu|+\lg V, Dv_\lambda+\nu\rg \qquad {\rm in }\; \R^{N+1}\; ,
$$
and
recall that there is some constant $\bar C$ independent of $\nu$ and $\lambda$ such that
$$
{\rm osc}( v_\lambda) \leq \bar C\qquad {\rm and }\qquad \|\lambda v_\lambda-\bar H(\nu)\|_\infty\leq \bar C\lambda\;.
$$

It is immediate that
$$
z(x,t)= v_\lambda(x,t)-v_\lambda(0)+\lg \nu, x\rg+(\bar H(\nu)+\bar C\lambda) t +C\;,
$$
with
$$
C= \bar C - \min_{y\in K_0} \lg \nu, x\rg +1\;,
$$
is a super-solution of (\ref{eq:geoWulff}). Since (\ref{eq:geoWulff}) geometric, it is immediate that $\max(z,0)$
is also a super-solution. Moreover, we clearly have
$\max(z(\cdot,0),0)\geq {\bf 1}_{K_0}$.

The standard comparison  gives, for all $(x,t)\in \R^N\times [0,+\infty)$,
$$
\max(z(x,t),0)\geq {\bf 1}_{K(t)}(x) \;,
$$
which, in turn,  implies that, for all $t\geq0$,
$$
K(t)\subset \{x\in \R^N\;:\; z(x,t)\geq 1\} \;.
$$

Recalling that $\bar H(\nu)\geq |\nu|=1$, we find, for all $(x,t)\in \R^N \times [1,\infty)$,
$$
z(x,t)\leq \lg \nu, x\rg+(\bar H(\nu)+\bar C\lambda) t +C+\bar C\leq \lg \nu, x\rg+\bar H(\nu)((1+\bar C\lambda)t +C+\bar C) \;,
$$
and, therefore, for all $t\geq0$,
$$
K(t) \subset \{x\in \R^N\; : \; \lg \nu, x\rg+\bar H(\nu)((1+\bar C\lambda)t +C+\bar C-1) \geq 0\} \qquad \forall t\geq 1\;.
$$

Letting  $\lambda\to 0$, we get, for all $t\geq0$ and a new positive constant $C$,
$$
K(t)\subset \{x\in \R^N\;:\; \lg \nu, x\rg+\bar H(\nu)(t +C) \geq 0\} \qquad \forall t\geq 1\;.
$$

Taking the intersection of the right-hand side over all $\nu$ we obtain, by the definition of ${\mathcal W}$, that, for all $t\geq0$,
$$
K(t)\subset (t+ C) {\mathcal W} \;.
$$

The proof of (\ref{InclWulff2}) is more intricate.
To this end, let $\bar R$, $r$ and $T_1$ be defined
by Lemma \ref{lem:keygrowth} and let $(K(t))_{t\geq0}$ be the solution of the front propagation problem starting at $K_0\subset \R^N$, a
compact set which contains
$Q_{\bar R}(k)$ for some $k\geq \Z^N$.
Then, for all $t\geq0$, 
\be\label{eq:inc}
B_{rt}\subset K(t+T_1).
\ee

Recall that  $u^\ep(x,t)={\bf 1}_{K(t/\ep)}(x/\ep)$
is the solution to
\be\label{eq:epWinclK}
\left\{\begin{array}{ll}
& u_t^\ep= |Du^\ep|+\lg V(x/\ep,t/\ep), Du^\ep\rg \quad \text{ in } \quad \R^N\times (0,\infty),\\
& u^\ep(x,0)= {\bf 1}_{K_0}(x/\ep).
\end{array}\right.
\ee
and, in view of \eqref{eq:inc}, for all $(x,t)\in\R^{N+1}$,
$$
u^\ep(x,t) \geq {\bf 1}_{B_{r(t/\ep-T_1)}}(x/\ep)= {\bf 1}_{B_{r(t-\ep T_1)}}(x)\;.
$$

Fix next $\delta\in(0,1)$ such that $T_1+\delta/(r\ep)$ is an integer and, for
$w_\delta(x)=(\delta-|x|) \vee 0$, let $w^\ep_\delta$ be the solution of (\ref{eq:epWinclK}) with initial
datum $w_\delta$. Then we know from Theorem \ref{theo:error} that there exists a constant $\bar C>0$ such that,
for all $t\in(0,1)$,
$$
\|w_\delta^\ep-\bar w_\delta\|_\infty \leq \bar C\ep^{1/3}
$$
where $\bar w_\delta$ is the solution of the homogenized problem

$$
\left\{\begin{array}{ll}
& \bar w_{\delta,t}= \bar H( D\bar w_\delta) \quad {\rm in} \quad \R^N\times (0,1)\;, \\
& \bar w_\delta(\cdot,0)= w_\delta \quad {\rm on} \quad \R^N \;.
\end{array}\right.
$$

Note that the constant $\bar C$ is independent of $\delta$ because the $w_\delta$'s have Lipschitz 
constant which are bounded uniformly in $\delta$.

From the Lax-Oleinik formula, $\bar w_\delta$ is given, for $(x,t)\in\R^N\times[0,\infty)$ by
$$
\bar w_\delta(x,t)= \sup_{y\in {\mathcal W}} w_\delta(x-ty)\;.
$$

Since $w_\delta\leq {\bf 1}_{B_{\delta}}\leq u^\ep(\cdot,\ep T_1+\delta/r)$, it follows, from the time-periodicity of $V$ and the choice of $\delta$,
that,  for all $(x,t)\in\R^N\times[0,1])$,
$$
w_\delta^\ep(x,t) \leq u^\ep(x,t+\ep T_1+\delta/r) \;.
$$

Hence, for all $t\in[0,1]$,
$$
\{w_\delta^\ep(\cdot,t)\geq \delta/2\}\subset \{u^\ep(\cdot,t+\ep T_1+\delta/r)\geq \delta/2\},
$$
while
$$
\{w_\delta^\ep(\cdot,t)\geq \delta/2\}  \supset \{\bar w_\delta(\cdot,t)\geq \delta/2+\bar C\ep^{1/3}\}
\supset \{\sup_{y\in {\mathcal W}} w_\delta(\cdot -ty) \geq \delta/2 +\bar C\ep^{1/3}\}\;.
$$

If we choose $\delta/2 -\bar C\ep^{1/3}>0$, then
$$
\{\sup_{y\in {\mathcal W}} w_\delta(\cdot -ty) \geq \delta/2 +\bar C\ep^{1/3}\} \supset
\{x\; :\; \sup_{y\in {\mathcal W}} (\delta-|x-ty|) \geq \delta/2 +\bar C\ep^{1/3}\} \supset t{\mathcal W}\;.
$$

Therefore, for all $t\in[0,1]$,
$$
\ep K((t+\ep T_1+\delta/r)/\ep) = \{u^\ep(\cdot,t+\ep T_1+\delta/r)\geq \delta/2\} \supset  t{\mathcal W} \;,
$$
i.e., for all $t\in[0,1/\ep]$,
\be\label{eq:inclincl}
K(t+T_1+\delta/(r\ep)) \supset t{\mathcal W}\;.
\ee

Finally, for $t$ sufficiently large, choose $\ep=1/(t-4\bar C t^{2/3}/r)$ and $\delta=(n-T_1)r\ep$
where $n$ is the integer part of $4\bar C t^{2/3}/r+1$. Then $n=T_1+\delta/(r\ep)$ is an integer, 
$\delta/2-\bar C \ep^{1/3}$ is positive and, applying inclusion (\ref{eq:inclincl}) to $t-n$ which belongs 
to $[0,1/\ep]$, we get 
$$
K(t) \supset (t-Ct^{2/3}) {\mathcal W}
$$
for some new constant $C$. 
\QED

We conclude the section with the proofs of the two lemmas used in the proof.

We have:

\noindent {\bf Proof of Lemma \ref{envahisseur} : } Fix $x\in\R^N$ and let $(K(t))_{t\geq 0}$ be the solution of the front propagation
$V_{x,t}= 1-V(x,t)$ starting from $K_0=\{x\}+\Z^N$, i.e., $K(t)=\{y\in \R^N\;:\; u(y,t)=0\}$ where $u$ is the solution to
$$
\left\{\begin{array}{l}
\partial_t u=|Du|+\lg V,Du\rg \quad \text{ in } \quad \R^N \times (0,\infty)\\
u= -d_{K_0} \quad \text{ on } \R^N ,
\end{array}\right.
$$
where $d_{K_0}$ is the distance function to the set $K_0$.

Then, for each $t>0$, the set $K(t)$ is $\Z^N$-periodic and a non empty interior (because it has an interior ball property, as
recalled in Appendix). So 
$\rho(t)=|K(t)\cap Q_1(x)|$ is positive for positive time and, following
the computation in the proof of Theorem \ref{theo:mainperiodic}, it satisfies for all $t_2>t_1\geq 0$,
$$
\rho(t_2)-\rho(t_1) \geq \frac{1}{c_I} \int_{t_1}^{t_2} \left(\min\{\rho(t), 1-\rho(t)\}\right)^{(N-1)/N}dt\; .
$$

Hence there exists a time $T$ depending only on $N$
such that $|K(T)\cap Q_1(x)|=1$. This means that $Q_1(x)\subset K(T)$.

It follows from the finite speed of propagation, that
there exits a positive integer $n_0$ such that the solution $\hat K(t)$ starting from $\{x+k\in \Z^N\;:\; |k|\leq n_0\}$
coincides with $K(t)$ on $Q_1(x)\times [0,T]$. Then $Q_1(x)\subset \hat K(T)$.
\QED

\noindent{\bf Proof of Lemma \ref {lem:keygrowth} : }  Consider
the solution $\bar u$  to
$$
\left\{\begin{array}{l}
\bar u_t=\bar H(D\bar u) \quad \text{ in } \quad \R^N \times (0,\infty),\\
\bar u(x,0)=(-|x|)\vee (-1)\quad \text{ on } \quad \R^N \times\{0\}\ .
\end{array}\right.
$$

Since, $\bar H(p)\geq |p|$,
given $0<\theta_1<\theta_0<1$ and $\delta>0$ small, there exists $\bar t\in(0,1)$ such that
$$
\{\bar u(\cdot,0)\geq -\theta_0\}+ B(0,\delta)\subset \{\bar u(\cdot,\bar t)\geq -\theta_1\} \;.
$$

The fact that the solution $u^\epsilon$ of

$$
\left\{\begin{array}{l}
\partial_t u^\epsilon=|Du^\epsilon|+\lg V(\frac{x}{\epsilon}, \frac{t}{\epsilon}), Du^\epsilon\rg \quad \text{ in } \quad \R^N\times(0,\infty),\\
u^\epsilon(x,0)=(-|x|)\vee (-1) \quad \text{ on } \quad \R^N \times\{0\}\; ,
\end{array}\right.
$$

\noindent converges, as $\ep\to0$, locally uniformly to $\bar u$, yields an   $\epsilon\in (0, \delta/4)$ such that
\begin{equation}\label{tutu}
\{u^\epsilon(\cdot,0)\geq -\theta_0\}+ B(0,\frac{\delta}{2})\subset \{u^\epsilon(\cdot,\bar t)\geq -\theta_0\}\;.
\end{equation}

Next fix  $n_0$ and $T$  as in Lemma \ref{envahisseur}, choose
$\epsilon$ such that $T_0= \frac{\bar t}{\epsilon}$ is an integer and
\be\label{epsmall2}
Q_{n_0+1}(0)\subset B(0, \frac{\delta}{2\ep})\;,
\ee
and set  $$\tilde K_0= \{x\ :\ u^\epsilon(\ep x,0)\geq -\theta_0\} \; .$$

The solution
of the front propagation problem starting from $\tilde K_0$ is given by $\tilde K(t)=\{x\;:\; u^\epsilon(\epsilon x,\epsilon t)\geq -\theta_0\}$.
From (\ref{tutu}) we have
$$
\tilde K_0 +B(0,\frac{\delta}{2\epsilon}) \subset \tilde K(T_0)\;,
$$
while from (\ref{epsmall2}), for any $k\in \Z^N$ with $|k|\leq n_0$, we have $\tilde K_0+k\subset \tilde K(T_0)$.

The periodicity of $V$ also implies that the solution of the front propagation problem
starting from $\tilde K_0+k$ is just $\tilde K(t)+k$ while the solution
starting from $\tilde K(T_0)$ is $K(t+T_0)$.

From the inclusion principle  we get, for all $k\in \Z^N$ with $|k|\leq n_0$ and all $t\geq 0$,
\begin{equation}\label{tyty}
\tilde K(t) +k \subset \tilde K(t+T_0) \;.
\end{equation}

Then Lemma \ref{envahisseur} implies that, for all $t\geq0$,
$$
\tilde K(t)+ Q_1 \subset \tilde K(t+T_0+T)\;.
$$

In particular, by induction, we get, for all positive integers $n$ and all $t\geq0$,
$$
\tilde K(t)+ Q_n\subset \tilde K(t+n(T_0+T)) \;.
$$

Choose a positive integer $M$ such that, for all $t\in [0,T_0+T]$,  $\tilde K(t)\subset Q_M$. Then, for all
positive integers $n$ such that $n\geq M$ and all $t\in [0,T_0 +T]$,

$$
Q_{n-M} \subset \tilde K(t)+Q_n \subset \tilde K(t+n(T_0+T))\;,
$$
and, hence, there exist $r>0$ and $T_1>0$ such that, for all $t\geq0$,
$$
Q_{rt} \subset \tilde K(t+T_1) \;.
$$

Finally  choose a positive integer $\bar R$ such that $\tilde K_0 \subset Q_{\bar R}$. Then, for any compact initial set
$K_0$ such that $Q_{\bar R}(k) \subset K_0$ for some $k\in \Z^N$, we have $\tilde K_0+k\subset K_0$. Therefore the solution of the front
propagation problem $K(t)$ starting from $K_0$ satisfies, for all $t\geq0$,
%\be\label{AAA}
$$
Q_{rt}(k) \subset \tilde K(t+T_1)+k \subset K(t+T_1) \;.
$$
%\ee
\QED

%%%%%%%%%%%%%%%%%%%%%%%%%%%%%%%%%%%%%%%%%%%%%%%%%%
%%%%%%%%%%%%%%%%%%%%%%%%%%%%%%%%%%%%%%%%%%%%%%%%%%
\section{Homogenization for $x$-dependent velocities at scale one }%$V=V(x,x/\ep,t/\ep)$} 
\label{sec:xtys}

Before we begin the proof, we remark that, since we are only able to prove that $\bar H$ is continuous with respect to the $x$-variable, uniqueness of the solution to
equation (\ref{eq:HomPbVxxep}) could be an issue. This is not, however, the case because, in view of (\ref{BoundIntV}) and (\ref{CroissVxxep}), $\bar H$
is coercive. Note that this is the reason we do not consider $V$
that also depends on a slow time variable, because then the coercivity of the averaged Hamiltonian would no longer
ensure a comparison principle for the limit problem. \\

We continue with the

{\bf Proof of Theorem \ref{HomVxxep}: } For any fixed $(x,P)\in \R^N\times \R^N$, let $v_\lambda^{P,x}=v_\lambda^{P,x}(y,t)$ be
the solution to
$$
v_{\lambda,t}^{P,x}+ \lambda v_\lambda^{P,x}=|Dv_\lambda^{P,x}+P|+\lg V(x,y,t), Dv_\lambda^{P,x}+P\rg\qquad \; {\rm in}\; \R^N\times \R\;.
$$

From Lemma \ref{osc} and Corollary \ref{cor:barH} we know that  there exist $\bar H(x,P)$ and, for any $M>0$,
a constant $C_M>0$ independent of $\lambda$ and $P$, such that, for all $x\in\R^N$ such that $|x|\leq M$,
\be\label{barHlambdavVxxep}
\|\bar H(x,P)- \lambda v_\lambda^{P,x}\|_\infty\leq C|P|\lambda \;.
\ee

Arguing as in the proof of Theorem \ref{theo:mainperiodic} one easily checks that
$\bar H$ is positively homogeneous of degree one and convex in $P$, and that
(\ref{CroissVxxep}) holds.

To complete the proof, it only remains to show that $\bar H$ is continuous in $(x,P)$.
To this end, observe that the
standard comparison arguments imply $\bar H$ is $(1+\|V\|_\infty)$-Lipschitz continuous in $P$.

Fix next $M>0$, $x_1,x_2\in \R^N$ with $|x_1|,|x_2|\leq M$. Once again the standard comparison of viscosity solutions (see \cite{b}) gives
\be\label{v1v2Vxxep}
\|v_\lambda^{P,x_1}-v_\lambda^{P,x_2}\|_\infty\leq C_M|P|\frac{|x_1-x_2|^{\lambda/L}}{\lambda}
\ee
where $L=\sup_{|x|\leq M} \|D_{y,s}V(x,\cdot,\cdot)\|_\infty$ and $C_M$ depends only on $M$.

Combining (\ref{barHlambdavVxxep}) with
(\ref{v1v2Vxxep}),
we obtain, for all $x_1,x_2 \in\R^N$ such that $|x_1|,|x_2|\leq M$ and $|x_1-x_2|\leq 1$,
$$
\begin{array}{l}
\ds{ |\bar H(x_1,P)-\bar H(x_2,P)|\; \leq  \; C_M|P|\omega(|x_1-x_2|) }
\end{array}
$$
with (the modulus) $\omega$ given, for $r\in(0,0]$, by
$$
\omega(r) = \inf_{\lambda \in (0,1]} (\lambda^{-1} r^{\lambda/L}+\lambda) \;.
$$

The proof of  the continuity of $\bar H$ is now complete. \QED

%%%%%%%%%%%%%%%%%%%%%%%%%%%%%%%%%%%%%%%%%%%%%%%%%%%
%%%%%%%%%%%%%%%%%%%%%%%%%%%%%%%%%%%%%%%%%%%%%%%%%%%
%\section{Discussion} \label{sec:future}
%
%
%Let us first point out that our homogenization result can be extended, with the same proof, to equations of the form
%$$
%\left\{\begin{array}{l}
%(i)\quad u_t^\ep=H(\frac{x}{\ep}, Du^\ep)+\lg V(\frac{x}{\ep}, \frac{t}{\ep}),Du^\ep \rg  \quad {\rm in } \quad \R^N\times (0,T)\\
%%\ee
%%with initial condition
%%\be\label{eq:ic}
%(ii) \quad u^\epsilon = u_0 \quad \text { on } \quad \R^N\times\{0\}\ .
%\end{array}\right.
%$$
%where $H: \R^N\times \R^N\to \R$ is Lipschitz continuous and $\Z^N-$periodic with respect to the first variable, 1-positively homogeneous,
%convex and such that 
%$$
%H(x,P)\geq \alpha |P|\qquad \forall (x,P)\in . 
%
%There remain several intriging points: for one thing we do not know how to handle homogenization of the $G-$equation in
%random, and even in almost-periodic, environment. 
%

%%%%%%%%%%%%%%%%%%%%%%%%%%%%%%%%%%%%%%%%%%%%%%%%%%
%%%%%%%%%%%%%%%%%%%%%%%%%%%%%%%%%%%%%%%%%%%%%%%%%%

\section{Appendix} %proof of Lemma \ref{EstiPeri}}

We present here the

\noindent{\bf Proof of Lemma \ref{EstiPeri}:}
Fix $\theta\in (\ \inf z(\cdot, 0), \sup z(\cdot, 0)\ )$ and let $K_{\theta}(t)=\{z(\cdot, t)\geq \theta\}$.
Since $z$ is a solution to
(\ref{evolz})-(i), $K_{\theta}(t)$ is given by
%representation formula for the set $K_\theta(t):=\{z(\cdot, t)\geq \theta\}$:
$$
K_{\theta}(t)=\left\{ x\in \R^N \ : \ \begin{array}{l} \exists \xi:[0,t]\to\R^N\ \mbox{\rm absolutely continuous such that }\xi(t)=x,\\
z(\xi(0),0)\geq \theta, \ |\xi'(s)+V(\xi(s),s)|\leq 1 \; {\rm a.e.}\ s\in (0,t)\end{array} \right\}\;,
$$
i.e., $K_{\theta}(t)$ is the reachable set for the controlled system
$$
\xi'(s)=f(\xi(s),s,\alpha(s)), \qquad |\alpha(s)|\leq 1 \; ,
$$
where $f(x, s, \alpha)=\alpha+V(x,s)$ and  $\alpha \in \R^N$ are such that $|\alpha|\leq 1$.

Since $V$ is of class ${\mathcal C}^{1,1}$, it follows from \cite{cf06} that,
for any $0<\tau<T$, there exits a constant $r=r(\tau, T)>0$ such that, for all $t\in [\tau, T]$, the set $K_\theta(t)$
has the interior ball property of radius $r$, i.e.,
$$
\text{ for all} \quad  x\in \partial K_\theta (t) \quad \text{ there exists } \quad y\in \R^N \quad  \mbox{\rm such that } \quad x\in \overline{B}(y,r) \subset K_\theta(t),
$$
where $\overline{B}(y,r)$ stands for the closed ball of radius $r$ centered at $y$, and, hence, for all $t>0$,
$K_{\theta}(t)$
is a set of finite perimeter (see, for instance, \cite{acm05}).\\
Note that \cite{cf06} deals
with time independent dynamics, but the proofs can be easily adapted to the time-dependent framework considered here.

Next set $d(x,t)=d_{K_{\theta}(t)}(x)$. It follows that $d$ is Lipschitz continuous in $(x,t)$ and, moreover,
\be\label{eq:d}
d_t = -1+\lg V(x-d Dd ,t)\ ,\  Dd(x,t)\rg \quad \mbox{\rm a.e in } \quad \{d>0\}\;.
\ee
 
It is, of course, clear  that $x\to d(x,t)$ is $1-$Lipschitz continuous. The $\|V\|_\infty+1$-Lipschitz continuity of
the map $t\to d(x,t)$  comes from the above representation formula of $K_{\theta,t}$.

To check (\ref{eq:d}), recall that, since the map $P\to |P|+\lg V,P\rg$ is convex,
$z$ is a sub-solution of (\ref{evolz})(i) with an equality (see \cite{bj90} and \cite{f89}), i.e., for any test function $\varphi$ and any local maximum point $(x,t)$ of $z-\varphi$, we have
$$
\varphi_t=|D\varphi|+\lg V, D\varphi\rg\;.
$$

The invariance property of the geometric equation (\ref{evolz})-(i) (see \cite{bss93}) then  implies that the map
$(x,t)\to {\bf 1}_{K_{\theta}(t)}(x)$ is also a sub-solution of (\ref{evolz})(i) with equality.

Assume now that $d$ is differentiable at some point $(\bar x,\bar t)$ and let $\bar y$ be the unique projection
of $\bar x$ onto $K_{\theta}(\bar t)$ and note that $w(y,t):={\bf 1}_{K_{\theta}(t)}(y)+d(y+\bar x-\bar y,t)$
has a local maximum at $(\bar y, \bar t)$. Indeed, if $y\notin K_{\theta}(t)$, then, for all  $(y,t)$ sufficiently close to $(\bar y, \bar t)$,
$$
w(y,t)= d(y+\bar x-\bar y,t) \leq 1+d(\bar x,\bar t)=w(\bar y,\bar t),
$$
and, if $y\in K_{\theta,t}$, then
$$
w(y,t)= 1+d(y+\bar x-\bar y,t) \leq 1+ |y+\bar x-\bar y-y|= 1+|\bar x-\bar y|=w(\bar y,\bar t)\;.
$$

Since $(x,t)\to {\bf 1}_{K_{\theta}(t)}(x)$ is a sub-solution of (\ref{evolz})(i) with equality and $|Dd(\bar x,\bar t)|=1$
we obtain
$$
-d_t(\bar x,\bar t)=1+\lg V(\bar y,\bar t)\ ,\ -Dd(\bar x,\bar t)\rg\;,
$$
and (\ref{eq:d}) holds because $\bar y=\bar x-d(\bar x,\bar t)Dd(\bar x,\bar t)$. \\

Next we  assume that $\theta$ is such that $\partial\{ z(\cdot, 0)\geq \theta\}=\{ z(\cdot, 0)= \theta\}$, which is  true
for almost all $\theta\in ( \inf z(\cdot, 0), \sup z(\cdot, 0))$.
It then follows from \cite{bss93} that, for all $t\geq0$,
$\partial\{ z(\cdot, t)\geq \theta\}=\{ z(\cdot, t)= \theta\}$ for any $t\geq0$.
\\

We now prove that, for all $\varphi\in C_c^{\infty}(\R^N)$  and all $t\in (0,T)$,
\be\label{appendix:klaim}
\int_{K_\theta(t)}\varphi(x)dx-\int_{K_\theta(0)}\varphi(x)dx
= \int_0^t \int_{\partial K_\theta(s)} \varphi(x)(1-\lg V(x,s),\nu(x,s)\rg ) d{\mathcal H}^{N-1}(x)ds\; .
\ee

To this end, for  $h>0$ small, let $\zeta_h:\R\to [0,1]$ be such that
$\zeta_h(\rho)= 1$ if $\rho\leq 0$,
$\zeta_h(\rho)= 1-\rho/h$ if $\rho\in [0,h]$ and $\zeta_h(\rho)= 0$ if $\rho\geq h$.

Multiplying (\ref{eq:d}) by $\varphi\zeta_h'(d)$ (which makes sense because the sets
$\{d(\cdot, s)=0\}=\partial K_\theta(t)$ and $\{d(\cdot,s)=h\}$ have a zero measure since
$K_\theta(t)$ and $\{d(\cdot,s)\leq h\}$ have the interior ball property) and integrating over $\R^N\times (0,t)$ gives
\be\label{zetah}
\begin{aligned}
& \int_{\R^N}\varphi(x)\zeta_h(d(x,t))dx-\int_{\R^N}\varphi(x)\zeta_h(d(x,0))dx \\
&   =
\int_0^t\frac{-1}{h}\int_{ \{ 0 < d(\cdot, s) < h\} }\varphi(x)( -1+\lg V(x-d(x,s)Dd(x,s),s),Dd(x,s) \rg ) \ dx\ ds\;.
\end{aligned}
\ee

Note that

$$
\lim_{h\to 0}\ [\int_{\R^N}\varphi(x)\zeta_h(d(x,t))dx-\int_{\R^N}\varphi(x)\zeta_h(d(x,0))dx]
\;=\;
\int_{K_\theta(t)}\varphi(x)dx-\int_{K_\theta(0)}\varphi(x)dx.
$$

We now concentrate on the right hand side of (\ref{zetah}). Since $|Dd|=1$ a.e. in $\{d>0\}$, the co-area formula implies
$$
\begin{aligned}
& \int_0^t\frac{-1}{h}\int_{ \{ 0 < d(\cdot, s) < h\} } \varphi(x)(-1+\lg V(x-d(x,s)Dd(x,s),s),Dd(x,s) \rg ) \ dx\ ds\\
& =\;\int_0^t\frac{1}{h}\int_0^h \int_{ \{ d(\cdot, s) =\sigma \} } \varphi(x)(1-\lg V(x-\sigma \nu^\sigma(x,s),s),\nu^\sigma(x,s) \rg) \
d{\mathcal H}^{N-1}(x)\
d\sigma ,\ ds
\end{aligned}
$$

\noindent where $\nu^\sigma(x,s)$ is the measure theoretic outward unit normal to the set $\{d(\cdot, s)<\sigma\}$,
which has finite perimeter since it satisfies the interior ball property with  radius $r+\sigma$.

In order to complete the proof of (\ref{appendix:klaim})
we just need to use the following

\begin{Lemma}\label{lem:append-limit} Let $E$ be a closed subset of $\R^N$ with the interior ball property of radius $r>0$. Then, for all 
compactly supported in $x$ $\Phi\in C(\R^N\times {\mathcal S}^{N-1})$,
$$
\lim_{\sigma\to 0} \int_{\{d_E(\cdot)=\sigma\}} \Phi(x, \nu^\sigma(x))d{\mathcal H}^{N-1}(x)=
\int_{\partial E} \Phi(x, \nu(x))d{\mathcal H}^{N-1}(x)\;,
$$
where $d_E(x)$ stands for the distance of $x$ to $E$ and $\nu^\sigma(x)$ (resp. $\nu(x)$) is the measure theoretic outward  unit normal
to $\{d_E(\cdot)<\sigma\}$ (resp. to $E$) at $x\in \partial E$.
\end{Lemma}

\noindent {\bf Proof: } Set $E_\sigma=\{d_E(\cdot)\leq \sigma\}$ and denote
by $\Pi_\sigma$ the projection of $\partial E_\sigma$ onto $\partial E$. It is known that $\Pi_\sigma$ is uniquely
defined for ${\mathcal H}^{N-1}-$a.e. point $x\in \partial E_\sigma$.  Let $\mu_\sigma=
{\mathcal H}^{N-1}\lfloor \partial E_\sigma$ and
$\bar \mu_\sigma=\Pi_\sigma \sharp \mu_\sigma$, $\mu_0=\bar \mu_0={\mathcal H}^{N-1}\lfloor \partial E$. \\

The first claim is that $\bar \mu_\sigma$ is absolutely continuous with respect to $\bar \mu_0$.
For this let us first recall that, since $E$ has the interior ball property of radius $r$, the map $\Pi_\sigma^{-1}$
is well-defined on $F_\sigma:=\{x\in \partial E\; ;\; \exists y\in \partial E_\sigma
\; {\rm with}\; |x-y|=\sigma\}$ and that $\Pi_\sigma^{-1}$ is Lipschitz continuous, with constant at most $(r+\sigma)/r$,
on $F_\sigma$ (see \cite{acm05}). So, if $Z$ is a Borel subset of $\partial E$, then
$$
\bar \mu_\sigma(Z)= {\mathcal H}^{N-1}(\Pi_\sigma^{-1}(Z))\leq {\rm Lip}(\Pi_\sigma^{-1}){\mathcal H}^{N-1}(Z)=\frac{r+\sigma}{r}\bar \mu_0(Z)\;.
$$

In particular,   $\bar \mu_\sigma$ is absolutely continuous with respect to $\bar \mu_0$ and, if
$\displaystyle{ f_\sigma= \frac{d\bar \mu_\sigma}{d\bar \mu_0}  }$, then $f_\sigma$ is bounded by $\frac{r+\sigma}{r}$
${\mathcal H}^{N-1}-$a.e. in $\partial E$. \\

Next we note that, for any
$\Phi\in C(\R^N\times {\mathcal S}^{N-1})$,
\be\label{intEsigma}
\int_{\partial E_\sigma}\Phi(x,\nu^\sigma(x))d{\mathcal H}^{N-1}(x)\; =\;
\int_{\partial E}\Phi(y+\sigma \nu(y),\nu(y))f_\sigma(y)d{\mathcal H}^{N-1}(y)\;.
\ee

Indeed, the definitions of $\bar \mu_\sigma$ and $f_\sigma$, give
$$
\begin{aligned}
& \int_{\partial E}\Phi(y+\sigma \nu(y),\nu(y))f_\sigma(y)d{\mathcal H}^{N-1}(y)
\; =\;
\int_{\partial E}\Phi(y+\sigma \nu(y),\nu(y))d(\Pi_\sigma \sharp  \mu_\sigma)(y)\\
&\qquad  =\; \int_{\partial E_\sigma}\Phi(\Pi_\sigma(x)+\sigma \nu(\Pi_\sigma(x)),\nu(\Pi_\sigma(x)))d \mu_\sigma(x),
\end{aligned}
$$
which implies (\ref{intEsigma}) because
$\nu(\Pi_\sigma(x))=\nu^\sigma(x)$ and $\Pi_\sigma(x)+\sigma \nu(\Pi_\sigma(x))=x$.\\

In view of (\ref{intEsigma}), to complete the proof, it only remains to show that the $(f_\sigma)$'s converge, as $\sigma\to 0$,  to $1$ in $L^1(\partial E, {\mathcal H}^{N-1})$.

Applying (\ref{intEsigma}) with $\Phi=1$ gives
$$
{\rm Per}(E_\sigma)=\int_{\partial E} f_\sigma(y)d{\mathcal H}^{N-1}(y)\;.
$$

Since $f_\sigma\leq \frac{r+\sigma}{r}$ ${\mathcal H}^{N-1}-$a.e. in $\partial E$, the lower-continuity of the perimeter implies that
$\lim_{\sigma\to 0}{\rm Per}(E_\sigma)= {\rm Per}(E)$. Using again the inequality
$f_\sigma\leq r^{-1}r+\sigma$ which holds ${\mathcal H}^{N-1}-$a.e. in $\partial E$, we obtain, in the limit $\sigma\to 0$,
$$
\int_{\partial E} |1-f_\sigma(y)|d{\mathcal H}^{N-1}(y)
 \leq \left|1-\frac{r+\sigma}{r}\right| {\rm Per}(\partial E) +\int_{\partial E} (\frac{r+\sigma}{r}-f_\sigma(y))d{\mathcal H}^{N-1}(y) \to0 .
$$

%Th completes the proof of the Lemma.

\QED

%%%%%%%%%%%%%%%%%%%%%%%%%%%%%%%%%%%%%%%%%%%%%%%%%%
%%%%%%%%%%%%%%%%%%%%%%%%%%%%%%%%%%%%%%%%%%%%%%%%%%
%%%%%%%%%%%%%%%%%%%%%%%%%%%%%%%%%%%%%%%%%%%%%%%%%%

\end{document}